\documentclass[12pt,draft]{article}

\usepackage{amsfonts,amsmath,amssymb,theorem}
\usepackage[english]{babel}

\DeclareMathOperator{\cp}{cap}
\DeclareMathOperator{\Lip}{Lip}
\DeclareMathOperator{\ACL}{ACL}

\DeclareMathOperator{\loc}{loc}

\DeclareMathOperator{\supp}{supp}
\DeclareMathOperator{\dist}{dist}
\DeclareMathOperator{\diam}{diam} 
\DeclareMathOperator{\esssup}{ess\,sup}

\begin{document}

\centerline{\bf Mappings with bounded $(P,Q)$-distortion on Carnot groups}
\vskip 0.5cm

\centerline{\bf A.~Ukhlov and S.~K.~Vodopyanov}
\vskip 0.5cm

\centerline{ABSTRACT.} 
\bigskip
{\small We study mappings with bounded
$(p,q)$-distortion associated to Sobolev spaces on
Carnot groups. Mappings of such type have applications to
the Sobolev type embedding theory and classification of
manifolds. For this class of mappings, we obtain estimates of
linear distortion, and a geometrical description. We prove also
Liouville type theorems and give some sufficient conditions for removability of sets.} 
\vskip 0.5cm

\centerline{\bf Introduction}

\vskip 0.5cm

Let $\Omega$ be an open set on Carnot group $\mathbb G$. We study Sobolev mappings $f : \Omega\to \mathbb G$ of the class $W^1_{q,\loc}(\Omega)$ under the following condition:
the local $p$-distortion, $p\geq q\geq 1$,
$$
K_p(x,f)=\inf\{ k(x): |D_H f|(x)\leq k(x)J(x,f)^{\frac{1}{p}}\}
$$
is integrable in the power $\kappa$, $1/\kappa = 1/q-1/p$.

The modern approach to the Sobolev mappings theory is based on relations between these mappings, Sobolev spaces theory and the nonlinear potential theory. In [32] Yu. G.
Reshetnyak proved that  a non-constant mapping $f:\Omega\to\mathbb
R^n$, $\Omega\subset\mathbb
R^n$, belonging to Sobolev space $W^1_{n,\loc}(\Omega)$ and possessing bounded distortion
$$
K(x)=\inf \{k(x) : |Df|^n(x)\leq k(x)J(x,f)\}\in L_\infty(\Omega),
$$
is continuous discrete and
open. Note, that the continuity of mappings with bounded
distortion follows only from finiteness of the distortion
[51] (see [24] for another proof of this property). The
necessity of  study of topological properties of Sobolev mappings arises
in the nonlinear elasticity theory. In some problems of this theory
the uniformly boundedness of distortion of a mapping is too
restrictive [2, 3]. Typically, we have only integrability of 
distortion for Sobolev mappings. Modern development of the theory of mappings with
integrable distortion shows, that we have no topological
properties without additional analytical assumptions (see, for example
[14--21,\,23,\,25,\,42]). Typically, it is required that the change of variable equality in the Lebesgue integral holds and that the connected components of the inverse image $f^{-1}(y)$ are compact for each point $y\in f(\Omega)$ (quasilight mappings). Therefore we assume some topological properties 
for mappings with integrable $p$-distortion.

An another approach to the geometric function theory is based on
investigation of the distortion function in the space of functions with
bounded mean oscillation. In this case it is assumed that the 
distortion is majorized by some $BMO$-function (see, for example
[1,\,35,\,36]).

The mappings with integrable distortion which are connected with
Sobolev spaces ($(p,q)$-quasiconformal mappings) are studied in
[38,\,39,\,46,\,55--58].

Mappings of such type have applications in the Sobolev
type embedding theory [8--10,\,58] and classification of
manifolds [37].

The aim of this paper is study basic analytical and
geometrical properties of mappings with bounded
$(p,q)$-distortion, the validity of Liouville type theorems
and description of the removable sets on Carnot groups in terms of capacity. 

We call a continuous mapping $f: \Omega\to\mathbb
G$, where $\Omega$ is an open set in Carnot group $\mathbb G$, the
{\it mapping with bounded $(p,q)$-distortion}, $1\leq q\leq
p<\infty$, if

1) $f$ is open and discrete;

2) mapping $f$ possesses Luzin condition $\mathcal N$ (the image of a set
of measure zero is a set of measure zero);

3) $f$ belongs to the Sobolev class $W^1_{1,\loc}(\Omega)$;

4) $J(x,f)\geq0$ a. e. and $J(x,f) \in L_{1,\loc}(\Omega)$;

5) the local $p$-distortion

$$
K_p(x,f)=\inf\{ k(x): |D_H f|(x)\leq k(x)J(x,f)^{\frac{1}{p}}\}
$$

\noindent belongs to $L_{\kappa}(\Omega)$, where the number
$\kappa$ is defined from the relation $1/\kappa = 1/q-1/p$
($\kappa = \infty$ if $q=p$).

The value $K_{p,q}(f;\Omega)=\|K_p(x,f)\mid L_{\kappa}(\Omega)\|$
is called the coefficient of distortion of the mapping $f$ in the
open set $\Omega$.

We note that in the case $\mathbb G=\mathbb R^n$, $n\geq 2$, and
$p=q=n$, this class of mappings coincide with the space mappings
with bounded distortion (quasiregular mappings) [32,\,34]. In
this case a mapping belongs to the Sobolev space
$W^1_{n,\loc}(\Omega)$ and topological properties follows only
from the boundedness of the distortion.

The main results of the article are the Liouville type theorem

\vskip 0.3 cm

{\bf Theorem~A.} Let $f:\mathbb G\to\mathbb G$ be a mapping with
bounded $(p,q)$-distortion, $\nu-1<q\leq p\leq\nu$. Then
$\cp\bigl(\mathbb G\setminus f(\mathbb G);W^1_s(\mathbb G)\bigr)=0$ where
$s={{p}/{(p-(\nu-1))}}$.

\vskip 0.3 cm

and the theorem about removable sets

\vskip 0.3 cm

{\bf Theorem~B.} Let $f:\Omega\setminus F\to\mathbb G$ be a
 mapping with bounded $(p,q)$-distortion, $p\geq
q\geq\nu$, and $F$ be a closed set in the domain $\Omega$,
$\cp\bigl(F;W^1_s(\mathbb G)\bigr)=0$, $s=p/(p-(\nu-1))$. Then

$1)$ the case $p\geq q>\nu$: the mapping $f$ extends to a continuous
mapping $\tilde f:\Omega\to\mathbb G$;

$2)$ the case  $p=q=\nu$: if $\cp\bigl(\mathbb Cf(\Omega\setminus
F);W^1_{\nu}(\mathbb G)\bigr)>0$ then the mapping $f$ extends to a
continuous mapping $\tilde f:\Omega\to\mathbb G\cup{\infty}$.
$($Hereafter  $\mathbb G\cup\{\infty\}$ is the one-point
compactification $\mathbb G$ with the standard topology.$)$

\vskip 0.3 cm

Note that removability sets for quasiregular mappings (the case $\mathbb G=\mathbb R^n$ and $p=q=n$) and 
mappings with bounded $(p,p)$-distortion (the case $n-1<p<n$)
was considered in [34] and [53] respectively.

The main technical tool of the article is 

\vskip 0.3 cm

{\bf Theorem~C.} Let $\Omega$ be an open set in $\mathbb G$ and
$f:\Omega\to \mathbb G$ be a mapping with bounded $(p,q)$-distortion.
The push-forward  function $v=f_{*}u: f(\Omega)\to \mathbb
R$, defined by the function $u\in C_0\cap W^1_p(\Omega)$ as
\begin{equation}
v(x)=
\begin{cases}
\Lambda \sum\limits_{z\in f^{-1}(x)} i(z,f)u(z),
\,\,\, & x\in f(\supp u),\\
0, \,\,\, & x\notin f(\supp u),
\end{cases}
\nonumber
\end{equation}
has following properties:

$1)$ $\supp v$ is a compact set and $f(\supp u)=\supp v$;

$2)$ $v$ is a continuous function;

$3)$ $v\in \ACL(f(\Omega))$ while $q>\nu-1$;

$4)$ in every compact embedded subdomain  $D\Subset\Omega$ the
inequality
$$
\|f_{\ast}u\mid L_{s}^{1}(f(D))\|\leq\Lambda
N(f,D)^{\frac{s-1}{s}} (K_{p,q}(f;D))^{\nu-1}\|u\mid
L_{r}^{1}(D)\|
$$
holds with $s=p/(p-(\nu-1))$, $r=q/(q-(\nu-1))>0$.

\vskip 0.3 cm

Recall that a stratified homogeneous group [7], or, in
other terminology, a Carnot group [31] is a~connected simply
connected nilpotent Lie group~ $\mathbb G$ whose Lie algebra~ $V$ is
decomposed into the direct sum~ $V_1\oplus\cdots\oplus V_m$ of
vector spaces such that $\dim V_1\geq 2$, $[V_1,\ V_k]=V_{k+1}$
for $1\leq k\leq m-1$ and $[V_1,\ V_m]=\{0\}$. Let
$X_{11},\dots,X_{1n_1}$ be left-invariant basis vector fields of
$V_1$. Since they generate $V$, for each $i$, $1<i\leq m$, one
can choose a basis $X_{ij}$ in $V_i$, $1\leq j\leq n_i=\dim V_i$,
consisting of commutators of order $i-1$ of fields $X_{1k}\in
V_1$. We identify elements $g$ of $\mathbb G$ with vectors $x\in\mathbb
R^N$, $N=\sum_{i=1}^{m}n_i$, $x=(x_{ij})$, $1\leq i\leq m$,
$1\leq j\leq n_i$ by means of exponential map $\exp(\sum
x_{ij}X_{ij})=g$. Dilations $\delta_t$ defined by the formula
$\delta_t x= (t^ix_{ij})_{1\leq i\leq m,\,1\leq j\leq n_j}$, are
automorphisms of $\mathbb G$ for each $t>0$. Lebesgue measure $dx$
on $\mathbb R^N$ is the bi-invariant Haar measure on~ $\mathbb G$
(which is generated by the Lebesgue measure by means of the
exponential map), and $d(\delta_t x)=t^{\nu}~dx$, where the number
$\nu=\sum_{i=1}^{m}in_i$is called the homogeneous dimension of the group~$\mathbb G$.
The Lebesgue measure $|E|$ of a measurable subset $E$ of $\mathbb G$
is $\int_E~dx$.

Euclidean space $\mathbb R^n$ with the standard structure is an
example of an abelian group: the vector fields $\partial/\partial
x_i$, $i=1,\dots,n$, have no non-trivial commutation relations
and form a basis in the corresponding Lie algebra. One example of
a non-abelian stratified group is the Heisenberg group $\mathbb
H^n$. Its Lie algebra has dimension $2n+1$ and its center is
$1$-dimensional. If $X_1,\dots,X_n,Y_1,\dots,Y_n,T$ is a basis in
the Heisenberg algebra, then the only non-trivial commutation
relations are $\bigl[X_i,Y_i\bigr]=-4T$, $i=1,\dots,n$; all other
vanish.

A homogeneous norm on the group $\mathbb G$ is a continuous function
$\rho:\mathbb G\to [0,\infty)$ that is $C^{\infty}$-smooth on $\mathbb
G\setminus\{0\}$ and has the following properties:

(a) $\rho(x)=\rho(x^{-1})$ and $\rho(\delta_t(x))=t\rho(x)$;

(b) $\rho(x)=0$ if and only if $x=0$;

(c) there exists a constant $c>0$ such that
$\rho(x_1\cdot x_2)\leq
c(\rho(x_1)+\rho(x_2))$ for all $x_1,x_2\in \mathbb G$.

The system of basis vectors $X_1,X_2,\dots,X_n$ of the space
$V_1$ (here and throughout we set $n_1=n$ and $X_{i1}=X_i$, where
$i=1,\dots,n$) satisfies H\"ormander's hypoellipticity condition.
The {\it Carnot--Carath\'eodory distance}~ $d(x,y)$ between two
points $x,y\in\mathbb G$ is defined as the greatest lower bound
of lengths of all horizontal curves with endpoints $x$ and
$y$, where the length is measured in the Riemannian metric with
respect to which the vector fields $X_1,\dots,X_n$ are
orthonormal and a~horizontal curve is thought of as a~piecewise
smooth path whose tangent vectors belong to~ $V_1$. It can be
shown that $d(x,y)$ is a~left-invariant finite metric with
respect to which the group of the automorphisms~ $\delta_t$ is
a~dilatation group with a coefficient $t$:
$d(\delta_tx,\delta_ty)=td(x,y)$. By definition, we put
$d(x)=d(0,x)$.

Let $\mathbb G$ be a Carnot group with one-parameter dilatation
group $\delta_t$, $t>0$, and homogeneous norm $\rho$, and let
$E$ be a measurable subset of $\mathbb G$. As usual, we denote by
$L_p(E)$, $p\in [1,\infty]$, the space of pth-power
integrable functions with the standard norm
$$
\|u\mid
L_p(E)\|=\biggl(\int\limits_{E}|u(x)|^p~dx\biggr)^{\frac{1}{p}},
$$
if $p\in [1,\infty)$, and $\|u\mid
L_{\infty}(E)\|=\esssup_{E}|u(x)|$ for $p=\infty$. We
denote by $L_{p,\loc}(E)$ the space of functions
$f:\Omega\to \mathbb R$ such that $f\in L_q(K)$ for each compact
subset $K$ of $E$.

Let $\Omega$ be an open set in $\mathbb G$. The Sobolev space $W^1_p(\Omega)$, $1\leq p\leq\infty$,
$(L^1_p(\Omega)$, $1\leq p\leq\infty$) consists of the functions
$u:\Omega\to\mathbb R$ locally integrable in $\Omega$, having a
weak derivatives $X_i u$ along the vector fields $X_i$,
$i=1,\dots,n$, and a finite (semi)norm
$$
\|u\mid W^1_p(\Omega)\|=\|u\mid L_p(\Omega)\|+\|\nabla_H u\mid
L_p(\Omega)\|\,\,\,\,\,(\|u\mid L^1_p(\Omega)\|=\|\nabla_H u\mid
L_p(\Omega)\|),
$$
where $\nabla_H u=(X_1u,\dots,X_nu)$ is the horizontal subgradient of $u$.
Recall that a locally integrable function $v_i:\Omega\to\mathbb R$
is called the weak derivative of $u$ along the vector
fields $X_i$, $i=1,\dots,n$, if
$$
\int\limits_{\Omega}v_i\varphi~dx =
-\int\limits_{\Omega}uX_i\varphi~dx
$$
for any test function $\varphi\in C_0^{\infty}(\Omega)$. If
$u\in W^1_p(U)$ for each bounded open set $U$ such that
$\overline{U}\subset\Omega$ then we say that $u$ belongs to the
class $W^1_{p,\loc}(\Omega)$.

We say that a function
$u:\Omega\to\mathbb R$ is absolutely continuous on lines ($u\in
\ACL(\Omega)$) if for each domain $U$ such that
$\overline{U}\subset\Omega$ and each foliation $\Gamma_i$ defined
by a left-invariant vector field $X_i$, $i=1,\dots,n$, $u$ is
absolutely continuous on $\gamma\cap U$ with respect to
one-dimensional Hausdorff measure for $d\gamma$-almost every curve
$\gamma\in\Gamma_i$. Recall that the measure $d\gamma$ on the
foliation $\Gamma_i$ equals  the inner product $i(X_i)$ of the
vector field $X_i$ and the bi-invariant volume $dx$; see, for
instance, [6, 54].

Let $(\mathbb X,r)$ be a complete metric space, $r$ be a metric in
$\mathbb X$, and $\Omega$ be an open set in Carnot group~$\mathbb G$.
We say that a mapping $f:\Omega\to\mathbb X$ is in the class
$W^1_{p,\loc}(\Omega;\mathbb X)$, if the following conditions fulfil:

(A) for each $z\in\mathbb X$ the function $[f]_z : x\in\Omega\mapsto
r(f(x),z)$ belongs to the class $W^1_{p,\loc}(\Omega)$;

(B) the family of functions $(\nabla_H[f]_z)_{z\in\mathbb X}$ has a
majorant in the class $L_{p,\loc}(\Omega)$, that is, there exists
a function $g\in L_{p,\loc}(\Omega)$ independent of $z$ such that
$|\nabla_H[f]_z(x)|\leq g(x)$ for almost all $x\in\Omega$ and all
$z\in\mathbb X$.

If $\mathbb X=\mathbb R$, then this definition of mappings with values
in $\mathbb R$ is equivalent to those given above. For $\mathbb G=\mathbb
R^n$ one obtains Reshetnyak's definition of Sobolev mappings [33]. If $\mathbb X = \widetilde{\mathbb G}$ is an another
stratified group with one-parameter dilation group
$\tilde{\delta}$, a homogeneous norm $\tilde{\rho}$, and so on,
then we obtain the definition of Sobolev mappings between two
groups and denote this class by symbol
$W^1_{p,\loc}(\Omega;\widetilde{\mathbb G})$. In this case it is
convenient  to use an equivalent description of Sobolev classes
[41,\,43] involving only the coordinate functions.

We say that a mappings $f:\Omega\to\widetilde{\mathbb G}$,
$\Omega\subset\mathbb G$, is in the Sobolev class
$HW^1_{p,\loc}(\Omega;\widetilde{\mathbb G})$, $1\leq p<\infty$, if

(1) $\tilde{\rho}(f(x))\in L_{p,\loc}(\Omega)$;

(2) the coordinate functions $f_{ij}$ belongs to $\ACL(\Omega)$
for all $i,j$ and $f_{1j}$ belongs to $W^1_{p,\loc}(\Omega)$ for
all $j=1,\dots,n_1$;

(3) the vector

$$
X_kf(x)=\sum\limits_{1\leq i\leq\tilde{m},\,1\leq
j\leq\tilde{n_i}} X_k
f_{ij}(x)\frac{\partial}{\partial\tilde{x}_{ij}}
$$
belongs to $\widetilde{V}_1$ for all $x\in\Omega$, $k=1,\dots,n$.

In the next proposition we list the equivalent descriptions of
Sobolev classes [41,\,43].

{\it Let $\Omega$ be a domain in a stratified group $\mathbb G$.
Then the following assertions hold.

$(1)$ A mapping $f:\Omega\to\widetilde{\mathbb G}$ belongs to
$W^1_{p,\loc}(\Omega;\widetilde{\mathbb G})$ if and only if $f$ can
be redefined on a set of measure zero to belong to
$HW^1_{p,\loc}(\Omega;\widetilde{\mathbb G})$.

$(2)$ The mapping $f:\Omega\to\widetilde{\mathbb G}$ belongs to
$W^1_{p,\loc}(\Omega;\widetilde{\mathbb G})$ if and only if for each
function $u\in\Lip(\widetilde{\mathbb G})$ the composition $u\circ f$
belongs to $W^1_{p,\loc}(\Omega)$ and $|\nabla_H(u\circ
f)|(x)\leq \Lip f\cdot g(x)$, where $g\in L_{p,\loc}(\Omega)$ is
independent of $f$.}

Since $X_i f(x)\in\widetilde{V}_1$ for almost all $x\in\Omega$ [31],
$i=1,\dots,n$, the linear mapping $D_H f(x)$ with matrix
$(X_if_j(x))$, $i,j=1,\ldots,n$, takes the horizontal subspace $V_1$ to
the horizontal subspace $\widetilde{V}_1$ and is called the formal
horizontal differential of the mapping $f$ at $x$. Let $|D_H
f(x)|$ be its norm
$$
|D_Hf(x)|=\sup\limits_{\xi\in
V_1,\,\rho(\xi)=1}\tilde{\rho}(D_Hf(x)(\xi)).
$$

Smooth mappings with differentials respecting the horizontal
structure are said to be contact. For this reason one could say
that mappings in the class $W^1_{p,\loc}(\Omega;\widetilde{\mathbb
G})$ are (weakly) contact. It is proved in [43,\,54] that a
formal horizontal differential $D_H:V_1\to\widetilde{V}_1$
induces a homomorphism $Df:V\to\widetilde{V}$ of the Lie algebras
$V$ and $\widetilde{V}$, which is called the formal differential.
Hence $|Df(x)|\leq C|D_Hf(x)|$ [43] for almost all $x\in\Omega$ with
constant $C$ depending only on the group and the homogeneous norm.
Here we set
$$
|Df(x)|=\sup\limits_{\xi\in V,\,
\rho(\xi)=1}\tilde{\rho}(Df(x)(\xi)).
$$
If $\mathbb G=\widetilde{\mathbb G}$, then the determinant of the matrix
$Df(x)$ is called the (formal) Jacobian of the mapping $f$, it is
denoted by $J(x,f)$.

A mapping with bounded $(p,q)$-distortion has a finite distortion:
$Df=0$ almost everywhere on the set $Z=\{x: J(x,f)=0\}$, and in
accordance with the lemma~1, which is formulated below, the
mapping $f$ belongs to the class $W^1_{q,\loc}(\Omega)$.
Therefore, in the case $\mathbb G=\mathbb R^n$ and $p=q=n\geq 2$, the
mapping $f$ is a mapping with bounded distortion:
$$
\frac{|Df(x)|^n}{J(x,f)}\leq K<+\infty \quad \text{a. e. in}\quad
\Omega.
$$

For the first time quasiconformal mappings in non-Riemannian
spaces were introduced by G. D. Mostow in 1973 [29]. In the
proof of the rigidity theorem he used quasiconformal
transformations of an ideal boundary of some symmetric space. M.
Gromov [11] shown, that a geometry of the such ideal boundary
is modeled by a nilpotent Lie group equipped with
Carnot-Cara\-the\-odo\-ry metric, which is non-Riemannian. These
works stimulated an interest to the study of quasiconformal
mappings on Carnot Lie groups and Carnot-Caratheodory spaces.

For  studying  quasiconformal mappings on Carnot groups, P. Pansu
introduced the notion of the differential ($\cal P
$-differential) on Carnot groups as a homomorphism of the
corresponding Lie algebra [31]. Using this notion A. Koranyi
and H. M. Reimann constructed foundations for the theory of
quasiconformal mappings on the Heisenberg groups [22]. The
problems of the theory of quasiconformal mappings on Carnot
groups were considered in [40,\,52].

The theory of mappings with bounded distortion on Carnot groups
was developed in articles [40,\,41,\,44,\,47,\,53] (see an alternative
approach in [4,\,5,\,12]). Quasimeromorphic mappings on
homogeneous groups and a corresponding value distribution theory
properties were studied in [26].

Note, that in papers [30,\,59,\,60] were constructed non-trivial examples of contact mappings and quasiconformal mappings on Carnot groups and classes of Carnot groups with an infinite-dimensional family of contact maps.

If $f$ is a homeomorphism with bounded distortion, then we obtain
an analytic definition of the quasiconformal mappings on
stratified group from [40]. In this paper it was proved the equivalence of the weakest analytical definition of quasiconformal mapping 
to the different descriptions
of the quasiconformality from papers [31] and [22].

While $\nu-1<q=p<\infty$, then the mappings under consideration
has bounded $p$-distortion and its properties on Riemannian
spaces (stratified groups) are studied in [37,\,48] ([53]).

Note, that if $q\geq \nu$, then Luzin condition $\mathcal N$ follows from
analytical assumptions [45,\,47].

Recall that the multiplicity function of a mappings $f$ on
the set $A\subset\Omega$ is called the value
$$
N(f,A)=\sup\limits_{x\in\mathbb G}N(x,f,A)=\sup\limits_{x\in\mathbb
G}\sharp\{f^{-1}(x)\cap A\}.
$$

\vskip 0.3 cm

{\bf Lemma~1.} Let $f : \Omega\to\mathbb G$ is a mapping with
bounded $(p,q)$-distortion. Fix a compact domain $D\Subset\Omega$
and an arbitrary bounded domain $D'$ such that $D'\supset
f{(\overline{D}})$. Then for every function $u\in
W^1_{\infty}(D')$ the composition $u\circ f$ belongs to
$W^1_q(D)$ and the inequality
$$
\|u\circ f\mid L^1_q(D)\|\leq
K_{p,q}(f;D)\biggl(\int\limits_{D'}|\nabla_H u(y)|^p
N(y,f,D)~dy\biggr)^{1/p}
$$
holds. Particularly, a mappings $f$ belongs to
$W^1_{q,\loc}(\Omega)$.

\vskip 0.3 cm

{\sc Proof.} In the case $\mathbb G=\mathbb R^n$ the lemma was proved
in [42]. Let function $u$ belongs to $W^1_{\infty}(D')$, $D'\subset\mathbb G$. Then the composition $u\circ f$ belongs to the class $\ACL(D)$ and has derivatives $X_i(u\circ f)$, $i=1,...,n$ along horizontal vector fields $X_i$ almost everywhere in $D$. Since $f$ has finite distortion, we have
\begin{multline}
\|u\circ f\mid L^1_q(D)\|\leq \biggl(\int\limits_D(|\nabla_H u|(f(x))|D_Hf|)^q(x)~dx\biggr)^{\frac{1}{q}}\\
=
\biggl(\int\limits_{D\setminus Z} |\nabla_H u|^q(f(x))|J(x,f)|^{\frac{q}{p}}\frac{|D_Hf|^q}{|J(x,f)|^{\frac{q}{p}}}~dx\biggr)^{\frac{1}{q}}.
\nonumber
\end{multline}
Using the H\"older inequality we derive
$$
\|u\circ f\mid L^1_q(D)\|
\leq
\biggl(\int\limits_D \biggl(\frac{|D_Hf|^p}{|J(x,f)|}\biggr)^{\frac{q}{p-q}}~dx\biggr)^{\frac{p-q}{pq}}
\biggl(\int\limits_D|\nabla_H u|^p(f(x))|J(x,f)|~dx\biggr)^{\frac{1}{p}}.
$$
Now applying the change of variable formula we obtain the required estimate of the norm.
\vskip 0.3 cm

In the next assertion we formulate some basic properties of 
mappings with bounded $(p,q)$-distortion.  One can consider these
properties as a generalization of the corresponding properties of
mappings with bounded distortion.

Let a  continuous mapping $f:\Omega\to\mathbb G$ be open and discrete. The set in points such that the mapping $f$
is not locally homeomorphic,
is called the branch set, and is denoted by $B_f$.

\vskip 0.3 cm

{\bf Proposition~1.} Let $f:\Omega\to\mathbb G$ is a mapping
with bounded $(p,q)$-distortion. Then

$1)$ if $q>\nu-1$ then $f$ is $\cal P$-differentiable almost
everywhere in $\Omega$;

$2)$ if $q>\nu-1$ then $|f(B_f)|=0$;

$3)$ if $x\in \Omega\setminus B_f$  then in some neighborhood
$W$ of  $x$, the restriction $f : W\to V=f(W)$ is a
homeomorphism on $W$, and in the case $q>\nu-1$ the inverse mapping
$g=f^{-1} : V\to W$ belongs to the class $W^1_{s,\loc}(V)$ where
$s=p/(p-\nu+1)$, and the norm $|D_Hg(y)|$ meets the
inequality
\begin{equation}
|D_Hg(y)|\leq K_p(x,f)^{\nu-1}J(y,g)^{1/s} 
\end{equation}
at almost all point $y=f(x)\in V$;

$4)$ if $1\leq q\leq p\leq \nu$ then $f$ has
$N^{-1}$-property: the condition $|A|=0$ implies
$|f^{-1}(A)|=0$;

$5)$ if $1\leq q\leq p\leq \nu$ then Jacobian $J(x,f)>0$ almost
everywhere in $\Omega$;

$6)$ if $\nu-1<q\leq p\leq\nu$ then measure $|B_f|$ of the branch
set  vanishes;

$7)$ for every measurable set $A\subset\Omega$ and every
measurable function $u$ a function $u\mapsto u(y)N(y,f,A)$ is
integrable in $\mathbb G$ if and only if a function $(u\circ
f)(x)J(x,f)$ is integrable on $A$, and the equality
$$
\int\limits_A(u\circ f)(x)|J(x,f)|~dx = \int\limits_{\mathbb
G}u(y)N(y,f,A)~dy
$$
holds;

if a compact domain $D\Subset\Omega$ is such that $|f(\partial
D)|=0$, and $u$ is a measurable function, then the function $u\mapsto
u(y)\mu(y,f,A)$ is integrable in $\mathbb G$ if and only if the
function $(u\circ f)(x)J(x,f)$ is integrable in $A$ and the
equality
$$
\int\limits_A(u\circ f)(x)J(x,f)~dx = \int\limits_{\mathbb
G}u(y)\mu(y,f,A)~dy
$$
holds. Here $\mu(y,f,D)$ is the topological degree of the mapping
$f:D\to\mathbb G$ at  $y\notin f(\partial D)$.

\vskip 0.3 cm

{\sc Proof.} Every continuous, open and discrete mapping is a
quasi-monotone mapping in the sense of work [40].  In this
work was proved that every quasi-monotone mapping, belonging to
$W^1_{q,\loc}(\Omega)$, $q>\nu-1$, is $\cal P$-differentiable
almost everywhere. Assertion~7 was proved in [54] (see also
[43]).

It is known that if $\cal P$-differential does not vanish at
a point $x\in\Omega$, then $x$ can not be a branch point of the
mapping $f$. Therefore $\cal P$-differential vanishes almost
everywhere on the set $B_f$. So, from the change of variable formula in the Lebesgue integral, it 
follows that measure of the image of the branch set equals
zero. Assertion~2 is proved.

Assertion~3 was proved in [55, Theorem~9]. We give here a
short proof of inequality (1) only. Since the inverse mapping $g$
belongs to Sobolev class $W^1_{1,\loc}(V)$, then, for almost all $x$,  the
derivatives $X_ig$, $i=1,\dots,n$, exist at the point $y=f(x)$,
and composition $D_H g(f(x))\circ D_Hf(x)$ is the identical map.
In order to prove (1), we note that
$$
|D_H g(y)\leq \frac{|D_Hf(x)|^{\nu-1}}{J(x,f)}.
$$
Hence we have
$$
|D_H g(y)|\leq \biggl(\frac{|D_Hf(x)|}{J(x,f)^{\frac{1}{p}}}\biggr)^{\nu-1}J(x,f)^{\frac{\nu-p+1}{p}}=K_p(x,f)^{\nu-1}J(y,g)^{\frac{1}{s}}.
$$

Assertion~4 can be proved by taking into account the methods of
the papers [41] where it was proved for mappings with
bounded distortion, [55] where it was proved in the case
$\nu-1<q\leq p\leq\nu$ and $f$ to be a homeomorphism, and
[42] where it was proved in the case $\mathbb G=\mathbb R^n$.

Assertion~5 follows from the previous statements and from the
change of variable formula.

For proving assertion~6, we note that as long as $J(x,f)>0$
almost everywhere in $\Omega$ then the branch set has measure zero.
\vskip 0.3 cm

Let $f:\Omega\to\mathbb G$ be a continuous open and discrete
mapping. A domain $D\subset\Omega$ is called normal, if
$f(\partial D)=\partial f(D)$. A normal neighborhood of a
point $x\in\Omega$ is called a neighborhood $U\subset\Omega$ of
the point $x$ such that

1) $U$ is a normal domain,

2) $U\cap f^{-1}(f(x))=\{x\}$.

We denote by symbol $\mathbb C D$ a supplement for the set $D$. A
connection between the multiplicity function, the degree of a
mapping, and the index of a mapping is formulated in

\vskip 0.3 cm

{\bf Lemma~2. ([27])} Let $f$ be a continuous open
discrete and sense-preserving mapping.

$1)$ If $U$ is an open domain such that
$\overline{U}\subset\Omega$, then for all points $x\in {\mathbb
C}f(\partial U)$ the inequality $N(x,f,U)\geq \mu(x,f,U)$ holds,
and for all points $x\in {\mathbb C}f(\partial U)\cup (U\cap B_f))$
the equality $N(x,f,U)= \mu(x,f,U)$ holds;

$2)$ if $U$ is a normal domain then for all points $x\in
f(U\setminus B_f)$ the equality $N(x,f,U)=N(f,U)$ holds;

$3)$ if $U$ is a normal domain then for all points $x\in f(U)$
the equality $\mu(x,f,U)=\sum\limits_{j=1}^k i(x_j,f)$ holds,
where $k=N(x,f,U)$ and $\{x_1, \dots,x_k\}=f^{-1}(x)\cap U$;

$4)$ if $U$ is a normal domain then for all points $x\in f(U)$
the equality $N(f,U)=\mu(x,f,U)=\mu(f,U)$ holds;

$5)$ for all points $x\in U$, $i(x,f)=N(f,U)$ if and only if $U$
is the normal neighborhood of the point $x$;

$6)$ $x\in B_f$ if and only if $i(x,f)\geq 2$.

\vskip 0.3 cm

In the present work we will use the following notation. Symbol
$U(x,f,s)$ denotes the $x$-component of the set
$f^{-1}(B(f(x),s))$. If $x\in\Omega$, $0<r<d(x,\partial\Omega)$
and $0<s<d(f(x),\partial f(\Omega))$, then
$$
l(x,r)=\inf\limits_{d(x,w)=r}d(f(w),f(x)),
$$
$$
L(x,r)=\sup\limits_{d(x,w)=r}d(f(w),f(x)).
$$
The value
$$
H(x,f)={\varlimsup\limits_{r\to 0}}\frac{L(x,r)}{l(x,r)}
$$
is called the linear distortion of the mapping $f$ at
$x\in\Omega$.

The following  lemma is true for a wide class of metric spaces.
The proof of this lemma in the Euclidean space [13,\,34] is
carried on a general situation with obvious modifications.

\vskip 0.3 cm

{\bf Lemma~3.} Let $f:\Omega\to\mathbb G$ be a continuous, open
and discrete mapping. Then for every point $x\in\Omega$ a number
$\sigma_x>0$ exists, such that for every $0<r<\sigma_x$, the next
properties are fulfilled:

$1)$ $U(x,r)$ is a normal neighborhood of the point $x$;

$2)$ $U(x,r)=U(x,\sigma_x)\cap f^{-1}(B(f(x),r))$;

$3)$ if $r<\sigma_x$, then $\partial U(x,r)=U(x,\sigma_x)\cap
f^{-1}(S(f(x),r))$;

$4)$ the sets $\mathbb C U(x,r)$ and $\mathbb C\overline{U}(x,r)$ are
connected;

$5)$ the domain $U(x,r)\setminus\overline{U}(x,t)$ is a ring for
every numbers $0<t<r\leq\sigma_x$.

\vskip 0.3 cm

\bigskip

\centerline{\bf 1.~Estimates of the linear distortion}

\bigskip

A well-ordered triple $(F_0,F_1;D)$ of nonempty sets, where $D$ is
an open set in $\mathbb G$, and $F_0$, $F_1$ are compact subsets of
$\overline{D}$, is called a condenser in the group $\mathbb G$.

The value
$$
\cp_p(E)=\cp_p(F_0,F_1;D)=\inf\int\limits_D |\nabla_Hv|^p~dx,
$$
where the infimum is taken over all nonnegative functions $v\in
C(F_0\cup F_1\cup D)\cap L^1_p(D)$, such that $v=0$ in a
neighborhood of the set $F_0$, and $v\geq 1$ in a neighborhood of
the set $F_1$, is called the $p$-capacity of the condenser
$E=(F_0,F_1;D)$. If $U\subset\mathbb G$ is an open set, and $C$ is a
compact in $U$, then a condenser $E=(\partial U, C; U)$  will
be denoted by $E=(U,C)$. Properties of $p$-capacity in the
geometry of vector fields satisfying H\"ormander hypoellipticity
condition,  can be found in [49,\,50].

\vskip 0.3 cm

{\bf Theorem~1.} Let $f:\Omega\to\mathbb G$ be a mapping with bounded $(p,q)$-distortion. Suppose that
$A\subset\Omega$ is a normal domain for which $N(f,A)<\infty$.
Then for a condenser $E=(F_0,F_1;A)$ the inequality
$$
\cp_q^{\frac{1}{q}}(F_0,F_1;A)\leq K_{p,q}(f,A)N(f,A)^{\frac{1}{p}}\cp_p^{\frac{1}{p}}(f(F_0),f(F_1);f(A))
$$
holds.

\vskip 0.3 cm

{\sc Proof.} Let function $u$ be an admissible function for condenser $(f(F_0),f(F_1);f(A))$. Then $u\circ f$ is an
admissible function for the condenser $(F_0,F_1;A)$. Applying the
change of variable formula (Proposition~1) we have
\begin{multline}
\cp_q^{\frac{1}{q}}(F_0,F_1;A)\leq\biggl(\int\limits_A|\nabla_H(u\circ
f)|^q~dx\biggr)^{\frac{1}{q}}\\
\leq\biggl(\int\limits_A|\nabla_H(u)|^q(f(x))|D_Hf|^q~dx\biggr)^{\frac{1}{q}}
\leq
K_{p,q}(f,A)\biggl(\int\limits_A|\nabla_H(u)|^p(f(x))J(x,f)~dx\biggr)^{\frac{1}{p}}\\
\leq K_{p,q}(f,A)N(f,A)^{\frac{1}{p}}\biggl(\int\limits_{f(A)}|\nabla_H(u)|^p(y)~dy\biggr)^{\frac{1}{p}}.
\nonumber
\end{multline}
Since $u$ is an arbitrary admissible function, then the theorem
is proved.

\vskip 0.3 cm

In the next theorem we give a geometric description of a mapping
with bounded $(p,q)$-distortion. This result generalizes the
corresponding results proved for $(p,q)$-quasi\-con\-for\-mal mappings
in [55], and for mappings with bounded $p$-distortion in
[53].

\vskip 0.3 cm

{\bf Theorem~2.} 
Let $f:\Omega\to\mathbb G$ be a mapping with
bounded $(p,q)$-distortion, $\nu-1<q\leq p<\infty$. Then there exists a number $r_0>0$  such that, for all $0<r<\lambda r<r_0$ where
$\lambda>1$ is a fix constant, the inequality
$$
\overline{\lim\limits_{r\to
0}}\frac{L(x,r)r^{(\nu-q)/q}}{|f(B(x,\lambda r))|^{\frac{1}{p}}}\bigg/ K_{p,q}(f,B(x,\lambda r))\leq c~i(x,f)^{\frac{1}{p}}
$$
holds at all points $x\in\Omega$.

\vskip 0.3 cm

{Proof.} We choose a point $x\in\Omega$ and let $\sigma_x>0$
be a number from Lemma~3. Suppose that a  number $t>0$ is  such
that $L(x,t)<\sigma_x$, and $r_0>0$ is such that the inclusion $U(x,f,s)\subset B(x,t)$ holds for all
$0<s\leq r_0$.

In the domain $\Omega$, we choose balls $B(x,r)\subset
B(x,\lambda r)\subset B(x,t)$, $\lambda>0$. Let  $y_1\in
f(S(x,\lambda r))$ be a point such that $d(f(x),y_1)=L(x,r)$.
Take a point  $y_2\in f(S(x,\lambda r))$ which is the most
remote from the point $y_1$, and a point $y_3\in f(S(x,\lambda r))$
 which is the least remote from the point $y_1$.
We will use the notations $d_0=d(y_2,y_1)$ and $d_1=d(y_2,y_3)$.
 In the domain $f(\Omega)$, consider the continuums
$$
F'_0=\mathbb CB(y_2,d_0)\cap f(B(x,\lambda r))\quad\text{and}\quad F'_1=\mathbb
CB(y_2,d_1)\cap f(B(x,\lambda r)).
$$
It is clear that the function
$$
\eta(y)=\frac{1}{cL(x,r)}\min(\dist(y,F'_0),cL(x,r))
$$
is admissible
for the condenser $(F'_0,F'_1;f(B(x,\lambda r)))$. Here a
constant $c$ is defined from the condition
$cL(x,r)=\dist(F'_0,F'_1)$. Note  that sets $F_0=B(x,\lambda
r)\cap f^{-1}(F'_0)$ and $F_1=B(x,\lambda r)\cap f^{-1}(F'_1)$
intersect spheres $S(x,t)$ where $r<t<\lambda r$. Indeed, if we
consider an arbitrary horizontal curve in $F'_i$, $i=0,1$,
connecting $f(S(x,r))$ and $f(S(x,\lambda r))$, then by Lemma~3
there exists a curve in $f^{-1}(F'_i)$, $i=0,1$, connecting the
sets $S(x,r)$ and $S(x,\lambda r)$. In this notation the
condenser $(F'_0,F'_1; f(B(x,\lambda r)))$ is the image of the
condenser $(F_0,F_1;B(x,\lambda r)\setminus B(x,r))$ under the
mapping $f$.   

In work [40,\,50]  the following estimate for Teihm\"uller
capacity of the condenser \linebreak $(F_0,F_1; B(x,\lambda r)\setminus
B(x,r))$ was obtained
$$
\cp_q (F_0,F_1; B(x,\lambda r)\setminus B(x,r))\geq c_1 r^{\nu-q}.
$$
Applying Theorem~1 we have
\begin{multline}
c_1 r^{\frac{\nu-q}{q}}\leq \cp_q^{\frac{1}{q}}(F_0,F_1;
B(x,\lambda r)\setminus B(x,r))\\
\leq K_{p,q}(f;B(x,\lambda r)) i(x,f)^{\frac{1}{p}}\cp_p^{\frac{1}{p}}(F'_0,F'_1; f(B(x,\lambda r)))\\
\leq K_{p,q}(f;B(x,\lambda r)) i(x,f)^{\frac{1}{p}}\frac{|f(B(x,\lambda r))|^{\frac{1}{p}}}{cL(x,r)}.
\nonumber
\end{multline}
Hence
$$
\frac{L(x,r)r^{\frac{\nu-q}{q}}}{|f(B(x,\lambda r))|^{\frac{1}{p}}}\leq c_2K_{p,q}(f;B(x,\lambda r)) i(x,f)^{\frac{1}{p}}.
$$
Passing to the limit as $r\to 0$ we obtain the desired
inequality. The theorem is proved.

\vskip 0.3 cm

{\bf Corollary~1.} 
If $p\geq q\geq \nu$, then we can take
$\lambda=1$. 

In the case $q>\nu$ , the statement can be proved alone the line of the proof
of the theorem with taking into account the positivity of the capacity of the two
one-point sets in $L^1_q(\mathbb G)$. In the case $p=q=\nu$ it is
follows from following proved in [53].

\vskip 0.3 cm

{\bf Lemma~4.} Let $f:\Omega\to\mathbb G$ be a nonconstant
mapping with bounded $\nu$-distortion and with the
distortion  coefficient $K(f)$. Then for every point $x\in\Omega$ the estimate
$$
H(x,f)\leq C(\nu, K_{\nu,\nu}(f,\Omega) i(x,f))<\infty
$$
holds, where constant $K_{\nu,\nu}(f,\Omega) i(x,f)$  depends
on the homogeneous dimension $\nu$ and the product
$K_{\nu,\nu}(f,\Omega) i(x,f)$ only.

\vskip 0.3 cm

\bigskip

\centerline{\bf 2.~Capacity estimates}

\bigskip

Let $f:\Omega\to\mathbb G$ be a continuous, open, discrete,  and
sense-preserving mapping.  Assume that $x$ belongs to $\Omega$. Consider a
horizontal curve $\beta: [a,b]\to\mathbb G$ such that
$\beta(a)=f(x)$. A curve $\alpha:\Delta_c\to\Omega$, where $c\leq
b$, and $\Delta_c=[a,c)$ or $\Delta_c=[a,b]$, is called the
lifting of the curve $\beta$ with the origin at the point $x$, if
$\alpha(a)=x$ and $f\circ\alpha=\beta|_{[a,c)}$. A curve $\alpha$
is called full (maximal) lifting of the curve $\beta$, if
the domain of $\alpha$ coincide with $[a,b]$.

The next assertion holds.

\vskip 0.3 cm

{\bf Lemma~5.} Suppose that $D$ is a normal domain for a
continuous open discrete and sense-preserving mapping $f$,
and $y\in f(D)$. Let $f^{-1}(y)\cap D=\{x_1,\dots,x_k\}$, where
$k=N(f,D)$, and every point is counted according to the index
$i(x,f)$. If $\beta:[a,b]\to f(D)$ is a horizontal curve,
$\beta(a)=y$, then

$1)$ there exists full liftings $\alpha_1,\dots,\alpha_k$ of the
path $\beta$ such that $\alpha_i$ begins at the point $x_i$,
$i=1,\dots,k$;

$2)$ $\#\{l:\alpha_l(t)=\alpha_j(t)\}=i(\alpha_j(t),f)$ for every
point $t\in [a,b]$ and every number $1\leq j\leq k$;

$3)$ $f^{-1}(\beta(t))\cap D=\{\alpha_1(t),\dots,\alpha_k(t)\}$ for
every point $t\in [a,b]$.

\vskip 0.3 cm

{\sc Proof} of  Lemma~5 repeats almost verbatim  the proof of the corresponding
assertion in the Euclidean space [27,\,34].

\vskip 0.3 cm

Topological properties of the function $v(y)$, defined in
equality (2), are described in the following lemma. Proof of this
lemma repeats the proof from [27, Lemma~5.4] with the
obvious modifications.

\vskip 0.3 cm

{\bf Lemma~6.} The function $v$ has the following properties:

$1)$ $\supp v$ is a compact set and $f(\supp u)=\supp v$;

$2)$ $v$ is a continuous function.

\vskip 0.3 cm

We will prove differential properties of the function $v$.

\vskip 0.3 cm

{\bf Lemma~7.} The function $v$ is an $\ACL$-function in the open
set $f(\Omega)$.

\vskip 0.3 cm

{\sc Proof.} It is enough to show that $v$ is an $ACL$-function in
some neighborhood of every point of $\supp v$. Fix a
point $x_0\in\supp v$, and let $f^{-1}(x_0)\cap\supp
u=\{q_1,\ldots,q_s\}$. By Lemma~5, there exists a number $r_0$:
$0<r_0<d(x_0,\partial (f(\Omega))$ such that normal
neighborhoods $U(q_i,f,r_0)$, $i=1,\ldots,s$, are disjoint. We
choose a number $r_1\leq r_0$ such that $B(x_0,r_1)\cap
f\bigl(\supp
u\setminus\bigcup\limits_{i=1}^{s}U(q_i,f,r_0)\bigr)=\emptyset$.
Then components $f^{-1}(B(x_0,r_1))$, which intersects
$\supp u$, are sets $U(q_i,f,r_1)$, $i=1,\ldots,s$. We put
$U=\bigcup\limits_{i=1}^{s}U(q_i,f,r_1)=\bigcup\limits_{i=1}^{s}U_i$
where $U_i=U(q_i,f,r_1)$.

Let $x_0=e$. The general case reduces to the previous one by applying
left transformations. Fix a horizontal vector field $X_{\tau}$, and
let $\cal Y$ be a fibration generated by this field. We choose a
cube $Q=S\beta_0$, where $\beta_0=\exp_s X_{\tau}$, $|s|\leq M$
and $S$ is a hyperplane transversal to $X_{\tau}$:
$$
S=\{(a;b)\mid x_{1\tau}=0,\ \ |a|\leq M,\ \ |b|\leq M\}
$$
(here $a=(x_{1j})$, $1\leq j\leq n_1$; $b=(x_{ij})$, $1<i\leq m$,
$1\leq j\leq n_i$, and $M$ is a number such that $Q\subset
B(e,r_1)$).

For every point $z\in S$, we denote by $\beta_z$ an element of
the horizontal fibration $z\beta_0$ with the beginning at the point $z$.
Thus, a cube $Q$ is the union of all such segments of integral
lines. We consider a tubular neighborhood of the fiber $\beta_z$
with the radius $r$:
$$
E(z,r)=\beta_zB(e,r)\bigcap Q= \Bigl(\bigcup\limits_{\tau\in
\beta_z}B(\tau,r)\Bigr)\bigcap Q.
$$

We recall that the mapping $\Phi$ defined on open subsets
from $D$ and taking nonnegative values is called a {\it finitely
quasiadditive} set function [57] if

1) for any point $x\in D$, exists $\delta$,
$0<\delta<\dist(x,\partial D)$, such that $0\leq
\Phi(B(x,\delta))<\infty$ (here and in what follows
$B(x,\delta)=\{y\in\mathbb R^n: |y-x|<\delta\}$);

2) for any  finite collection $U_i\subset U\subset D$,
$i=1,\dots,k$, of mutually disjoint open sets the following
inequality $\sum\limits_{i=1}^k \Phi(U_i)\leq \Phi(U)$ takes
place.

Obviously, the inequality in the second condition of this
definition can be extended to a countable collection of mutually
disjoint open sets from $D$, so a finitely quasiadditive set
function is also {\it countable quasiadditive.}

If instead of second condition we suppose that for any finite
collection $U_i\subset D$, $i=1,\dots,k$, of mutually disjoint
open sets the equality
$$
\sum\limits_{i=1}^k \Phi(U_i)= \Phi(U)
$$
takes place, then such a function is said to be {\it finitely
additive}. If the equality in this condition can be extended to a
countable collection of mutually disjoint open sets from $D$,
then such a function is said to be {\it countable additive.}

A mapping $\Phi$ defined on open subsets form $D$ and taking
nonnegative values is called a {\it monotone} set function [57] if
$\Phi(U_1)\leq\Phi(U_2)$ under the condition that $U_1\subset
U_2\subset D$ are open sets.

Let us formulate a result from [57] in a form convenient for us.

\vskip 0.3cm

{\bf Proposition~2. ([57])} Let a finitely quasiadditive set
function $\Phi$ be defined on open subsets of the domain
$D\subset\mathbb R^n$. Then for almost all points $x\in D$ the
finite derivative
$$
\Phi'(x)=\lim\limits_{\delta\to 0, B_{\delta}\ni x}
\frac{\Phi(B_{\delta})}{|B_{\delta}|}
$$
exists and for any open set $U\subset D$, the inequality
$$
\int\limits_{U}\Phi'(x)~dx\leq \Phi(U)
$$
is valid.
\vskip 0.3cm

A nonnegative function $\Phi$ defined on a certain collection of
measurable sets from the open set $D$ and taking finite values is
said to be {\it absolutely continuous} if, for every number
$\varepsilon>0$, a number $\delta>0$ can be found such that
$\Phi(A)<\varepsilon$ for any measurable sets $A\subset D$ from
the domain of definition, which satisfies the condition
$|A|<\delta$.

Let the set function $\Phi$ be defined by the rule $\Phi(S)=K_{pq}(f;S)^{pq/(p-q)}$, $S\subset D^{\prime}$.
We define also a Borel function $\Psi$  as $\Psi(S)=|U\cap
f^{-1}(S\cap Q)|$. 
Fix a point $z\in S$ such that the upper volume
derivatives $\overline{\Psi'}(z)$ and $\overline{\Psi'}(z)$ are finite. It is sufficiently to prove that
the function $v$ is absolutely continuous on $\beta_z$. We
consider a lifting $\alpha:[a,b]\to U$ of the curve
$\beta_z:[a,b]\to Q$ such, that if $t_0\in [a,b]$ and $x_0\in U$,
$f(x_0)=\beta(t_0)$, then $\alpha(t_0)=x_0$ and
$f\circ\alpha=\beta$. Such lifting exists by Lemma~5. For
verifying the absolute continuity of  $v$ on $\beta_z$ we need
the following assertion. For proving this assertion, we apply
methods of paper [55].

We recall, that the Hausdorff $\alpha$-measure of a set $A$ is
said to be the value
$$
{\cal H}^{\alpha}(A)=\lim_{\varepsilon\to 0} \Bigl\{\inf\sum r_i^{\alpha}\Bigr\}
$$
where the infimum is taken over all coverings of the set $A$ by balls
$B_i$ with radii $r_i$ less then $\varepsilon$.

It is known [52, Proposition~1] that additive function
$\Psi$, defined on Borel subsets of $\mathbb G$, possesses the next
property: the upper volume derivative
$$
\overline{\Psi'}(z)=\varlimsup\limits_{r\to
0}\frac{\Psi(E(z,r))}{r^{\nu-1}}<+\infty
$$
exists for almost all $z\in S$.

\vskip 0.3 cm

{\bf Lemma~8.} Fix a point $z\in S$  in  which the functions 
$\Psi$ and $\Phi$ have a finite upper derivative
$$
\varlimsup\limits_{r\to 0}\frac{\Psi(E(z,r))}{r^{\nu-1}}<\infty, \quad \varlimsup\limits_{r\to 0}\frac{\Phi(E(z,r))}{r^{\nu-1}}<\infty.
$$
Let $\alpha(t):[a,b]\to U_k$ be a lifting of a horizontal curve
$\beta_z(t):[a,b]\to Q$ $($here $U_k\subset U$ is a fix normal
neighborhood$)$. Then $\alpha$ is an absolutely continuous curve
$($with respect to the Hausdorff $1$-measure$)$.

\vskip 0.3 cm

{\it Proof of the Lemma}~8 uses the following result.

\vskip 0.3 cm

{\bf Lemma~9 \rm[55, Proposition~5].} Let $E$ be a connected
set, and $G= \{x:d(x,E)\leq c_0\diam E\}$ where $c_0$ is a
small number depending on the constant in the generalized triangle
inequality. Then
$$
\cp_p^{\nu-1}(E,G)\geq c\frac{(\diam E)^p}{|G|^{p-(\nu-1)}}
$$
for $\nu-1<p<\infty$ where a constant $c$ depends only on $\nu$
and $p$.

\vskip 0.3 cm

On a horizontal curve $\beta_z$ we choose mutually disjoint closed
arcs
$[\delta_1,\overline\delta_1],\ldots,[\delta_l,\overline\delta_l]$
with lengths
$\Delta_1,\ldots,\Delta_l$ respectively such
that
$$
\sum\limits_{i=1}^{l}\Delta_i<\delta.
$$
By symbol $R_i$, we denote  the union of balls
$B_c(\beta_z(\tau),r))$ where
$\beta_z(\tau)\in[\delta_i,\overline\delta_i]$ (balls are
considered in Carnot-Charact\'erology metric $d_c$).
In this case the
set $(R_i,[\delta_i,\overline\delta_i])$ is a condenser. We choose
a small number $r>0$ such that the following properties fulfilled:
for some constant $c_1$ we have $r<c_1\Delta_i$, sets $R_i$, $i=1,\ldots,l$, are
disjoint, and the conditions of Lemma~9 hold.

Let
$[a_i,b_i]=\beta_z^{-1}([\delta_i,\overline\delta_i])\subset[a,b]$.
Then $\alpha([a_i,b_i])\subset f^{-1}(R_i)\cap U_k$ and the pair
$(f^{-1}(R_i)\cap U_k,\alpha([a_i,b_i]))$ is also a condenser,
since
$$
f^{-1}(R_i)\cap
U_k=\bigcup\limits_{\tau\in[a_i,b_i]}U(\alpha(\tau),f,r)
$$
is an open connected set. We note, that the image of the condenser
$E=(f^{-1}(R_i)\cap U_k,\alpha([a_i,b_i]))$ is the condenser
$f(E)=(R_i,[\delta_i,\overline\delta_i])$, since
$f\circ(\alpha([a_i,b_i]))=[\delta_i,\overline\delta_i]$ and
$$
f(f^{-1}(R_i)\cap U_k)=f
\Bigl(\bigcup\limits_{\tau\in[a_i,b_i]}U(\alpha(\tau),f,r)\Bigr) =
\bigcup\limits_{\tau\in[a_i,b_i]}B_c(\beta_z(\tau),r)=R_i.
$$

We note that the function $$u(q)=\frac{d(q,\partial R_i)}{r}$$ is
an admissible function for the condenser
$(R_i,[\delta_i,\overline\delta_i])$ and $|\nabla_{H}u(q)|\leq
1/r$. Hence, we obtain the inequalities
\begin{equation}
\cp_{p}(R_i,[\delta_i,\overline\delta_i]) \leq
\int\limits_{R_i}|\nabla_{H}u(x)|^{p}\,dx \leq
\frac{|R_i|}{r^{p}} \leq \frac{c_1\Delta_ir^{\nu-1}}{r^{p}}.
\end{equation}
On the other hand, by Lemma~9, we have
\begin{equation}
\cp_{q}(f^{-1}(R_i)\cap U_k,\alpha([a_i,b_i])) \geq
c_2\frac{\diam^{\frac{p}{\nu-1}}(\alpha([a_i,b_i]))}
{|f^{-1}(R_i)\cap U_k|^{\frac{q-\nu+1}{\nu-1}}}.
\end{equation}
Using Theorem~1 and the inequalities (3), (4) we have
$$
c_2^{\frac{1}{q}}\frac{\diam^{\frac{1}{\nu-1}}(\alpha([a_i,b_i]))}
{|f^{-1}(R_i)\cap U_k|^{\frac{q-\nu+1}{q(\nu-1)}}} \leq
K_{pq}(f;R_i)N(f,R_i)^{\frac{1}{p}}c_1^{\frac{1}{p}}{\Delta_i}^{\frac{1}{p}} r^{\frac{\nu-1-p}{p}}.
$$
Hence,
$$
\diam(\alpha([a_i,b_i]))\leq c_3
\biggl(\frac{\Phi(R_i)}{r^{\nu-1}}\biggr)^{\frac{(p-q)(\nu-1)}{pq}}
\biggl(\frac{|f^{-1}(R_i)\cap
U_k|}{r^{\nu-1}}\biggr)^{\frac{q-\nu+1}{q}}
\Delta_i^{\frac{\nu-1}{p}}
$$
where $c_3=(c_1
N(f,R_i))^{\frac{\nu-1}{p}}/c_2^{\frac{\nu-1}{q}}$.

We note that a tubular
neighborhood
$$
E(z,\lambda^{\prime}r)=\{x\in Q\mid d(x,\beta_z)<
\lambda^{\prime}r\}
$$
where $\lambda^{\prime}$ is a some constant, depending on the
geometry of the group $\mathbb G$ only, contains the union
$\bigcup\limits_{i=1}^{l}R_i$. So, we have
$\bigcup\limits_{i=1}^{l}f^{-1}(R_i)\cap U_k\subset
f^{-1}(E(z,\lambda^{\prime}r))$. Summing the last inequality over
$i=1,\ldots,l$ and applying H\"older inequality, we obtain
$$
\sum\limits_{i=1}^{l}\diam(\alpha([a_i,b_i]))\leq c_4
\biggl(\frac{\Psi(E(z,\lambda^{\prime}r))}{(\lambda^{\prime}r)^{\nu-1}}\biggr)^{\frac{(p-q)(\nu-1)}{pq}}
\biggl(\frac{|f^{-1}(E(z,\lambda^{\prime}r)|}{(\lambda^{\prime}r)^{\nu-1}}\biggr)^{\frac{q-\nu+1}{q}}
\biggl(\sum\limits_{i=1}^{l}\Delta_i\biggr)^{\frac{\nu-1}{p}}.
$$

Letting $r$ go to $0$, and using the condition of Lemma at the
point $z$, we have
$$
\sum\limits_{i=1}^{l}\diam(\alpha([a_i,b_i])) \leq
c_5\biggl(\sum\limits_{i=1}^{l}\Delta_i\biggr)^{\frac{\nu-1}{p}}.
$$
Thus, the curve $\alpha$ is absolutely continuous. Lemma~8 follows.

For proving Lemma~7 it is sufficiently to verify that
the function $v$ is absolutely continuous on $\beta_z$, where
$z\in S$ and $\beta_z$ are the same as in Lemma~8. Let
$U_i=U(q_i,f,r_1)$ be normal neighborhoods defined at the
beginning of the proof Lemma~7, and $\beta=\beta_z$.
By Lemma~2, the relation
$$ \sum\limits_{x\in f^{-1}(\beta(t))\cap
U_i}i(x,f)= N(f,U_i)=i(q_i,f)=k(i)
$$
holds. For every $i=1,\ldots,s$, we chose full liftings
$\alpha_{i,j}$, $j=1,\ldots,k(i)$, of the curve $\beta$ in $U_i$
according to Lemma~5. By Lemma~8, all curves
$\alpha_{i,j}$ are absolutely continuous. For the fix number $i$,
$i=1,\ldots,k(i)$, we know number of liftings passing over a point
belonging to $f^{-1}(\beta(t))\cap U_i$: Lemma~5 implies that
$$
\sum\limits_{x\in f^{-1}(\beta(t))\cap U_i}
i(x,f)u(x)=\sum\limits_{j=1}^{k(i)}u(\alpha_{i,j}(t)).
$$
From here we come to
\begin{equation}
v(\beta(t))=\sum\limits_{i=1}^{s} \Lambda\sum\limits_{x\in
f^{-1}(\beta(t))\cap U_i} i(x,f)u(x)=
\Lambda\sum\limits_{i=1}^{s}\sum\limits_{j=1}^{k(i)}u(\alpha_{i,j}(t)).
\end{equation}
According to (5), it is sufficiently to prove that
$u(\alpha_{i,j}(t))$ is absolutely continuous on $[a,b]$ for
every $i,j$. The last property follows from the absolute
continuity of the curve $\alpha_{i,j}$ and from
Lipschitz continuity of the function $u$ in the domain $U$:
$$|u(z_1)-u(z_2)|\leq Ld(z_1,z_2),\quad  z_1,z_2\in U,$$ since $u\in
C^{1}_0$. Lemma~7 is proved.

\vskip 0.3 cm

{\bf Lemma~10.} Suppose that $x_0\in\supp v\setminus f(\supp
u\cap B_f)$. Then there exists a neighborhood $V_0$ of $x_0$
such that, for every connected neighborhood $V \subset V_0$ of
the point $x_0$, the following conditions holds:

$1)$ $V\cap f(\supp u\cap B_f)=\emptyset$;

$2)$ the number of components $f^{-1}(V)$ intersecting $\supp u$
is finite; we denote them by $D_1,\ldots,D_k$;

$3)$ the restriction $f\vert_{D_i}=f_i:D_i\to V$, $i=1,\ldots,k$,
is a $(p,q)$-quasiconformal homeomorphism;

$4)$ if $g_i=f^{-1}_i$, then $|\nabla_{H}
v(z)|\leq\Lambda\sum\limits_{i=1}^{k}|\nabla_{H} u(g_i(z))|
|D_Hg_i(z)|$ for almost all points $z\in V$.

\vskip 0.3 cm

{\sc Proof.} We chose neighborhoods $U_1,U_2,\ldots,U_k$
of points $f^{-1}(x_0)\cap\supp u$ such that $\overline
U_i\Subset\Omega$ and the restriction $f\vert_{\overline U_i}$ is
an injective mapping. We have to show that
$$
V_0=\biggl(\bigcap\limits_{i=1}^{k}f(U_i)\biggr) \setminus
f\biggl(\supp u\setminus\bigcup\limits_{i=1}^{k}U_i\biggr)
$$
is a desired neighborhood $V_0$ of the point $x_0$.

Let $V\subset V_0$ be a connected neighborhood of the point
$x_0$. The first assertion is valid since $U_i$ does not intersect
$B_f$ for every $i=1,\ldots,k$. If a connected component $D$ of the preimage
$f^{-1}(V)$ intersects $\supp u$ then it intersects
one of the neighborhoods $U_i$. Since the restriction
$f\vert_{\overline U_i}$ is injective then $V_0\cap f(\partial
U_i)=\emptyset$, and hence $D\cap\partial U_i=\emptyset$. From
here it follows  $D\subset U_i$. Thus, the second and the third
assertions of  Lemma~10 are proved.

Since the mappings $g_i=f^{-1}_i$ are $(r,s)$-quasiconformal
homeomorphisms [55], then they are $\cal P$-differentiable
almost everywhere in $V$. In view of  $i(q,f)=1$ at the points
$q\in\Omega\setminus B_f$ we have
$$
v(z)=\Lambda\sum\limits_{i=1}^{k}u(g_i(z))
$$
 for every points $z\in V$.
Then at every point $z\in V$ of the $\cal
P$-differentiability of all $g_i$, we have
$$
|\nabla_{H}v(z)|\leq\Lambda\sum\limits_{i=1}^{k}|\nabla_{H}
u(g_i(z))| |D_Hg_i(z)|.
$$
Hence, the last assertion of Lemma is also proved.

\vskip 0.3 cm

{\sl Proof of Theorem~3.} Let $u\in C_{0}^{1}(\Omega)$ be an
arbitrary function. Applying Vitali covering theorem, we
obtain a countable
collection of disjoint balls $\{B_1,B_2,\ldots\}$ ,
covering $f(D)\setminus f(\supp u\cap B_f)$ up to a set of measure zero, such that for balls
$B_j$ intersecting $\supp v$ the conditions 1--4 of the previous
Lemma hold. Since in view of Proposition~1 $|f(B_f)|=0$,
we have
$$
\int\limits_{f(D)}|\nabla_{H}v|^s\,dz \leq
\sum\limits_{j=1}^{\infty}\int\limits_{B_j}|\nabla_{H} v|^s\,dx,
$$
where $s=\frac{p}{p-(\nu-1)}$. For some index $j$ we fix ball $B_j$. If the ball
$B_j$ does not intersects $\supp v$, then
$\int\limits_{B_j}|\nabla_{H} v|^s\,dx=0$. If $B_j$ intersects
$\supp v$, then $g_i :B_j\to W_i$, $i=1,\ldots,k(j)$, is the
inverse mapping, defined by property~3 of Lemma~10. From
the Minkowski inequality and from the relation (1) it follows
\begin{multline}
\biggl( \int\limits_{B_j}|\nabla_{H} v|^s\,dz
\biggr)^{1/s}\\
\leq \biggr( \int\limits_{B_j}
\Bigl|\Lambda\sum\limits_{i=1}^{l}|\nabla_{H}
u(g_i(z))||D_hg_i(z)|\Bigl|^s\,dz \biggl)^{1/s}
\leq \Lambda\sum\limits_{i=1}^{l} \biggr(
\int\limits_{B_j}|\nabla_{H} u(g_i(z))|^s|D_hg_i(z)|^s\,dz
\biggl)^{1/s}
\\
\leq \Lambda\sum\limits_{i=1}^{l} \biggr(
\int\limits_{B_j}|\nabla_{H} u(g_i(z))|^sK_p(g_i(z);f)^{s(\nu-1)}
J(z,g_i)\,dz
\biggl)^{1/s}\\
\leq 
\begin{cases} \Lambda K_{p,p}^{\nu-1}(f;D)\sum\limits_{i=1}^{l}
\biggr( \int\limits_{g_i(B_j)}|\nabla_{H} u|^s\,dx
\biggl)^{1/s},\quad&\text{при}\ q=p,\\
\Lambda \sum\limits_{i=1}^{l}
\biggl(\int\limits_{g_i(B_j)}K_p(x;f)^{\frac{pq}{p-q}}\,dx
\biggr)^{\frac{r-s}{rs}} \biggr(\int\limits_{g_i(B_j)}|\nabla_{H}
u|^r\,dx \biggl)^{1/r},\quad&\text{при}\ q<p.
\end{cases}
\nonumber
\end{multline}
In the last inequality we have used the $(r,s)$-quasiconformality of
the mappings $g_i$. We can also assume that $\diam W_i<d(\supp
u,\partial D)$. Then $l\leq N(f,D)$ and, in the case $p=q$,
H\"older inequality implies
\begin{multline}
\biggr( \int\limits_{B_j}|\nabla_{H} v|^s\,dz
\biggl)^{1/s} \leq \Lambda K^{\nu-1}_{p,p}(f;D)l^{(s-1)/s} \biggr(
\sum\limits_{i=1}^{l}\int\limits_{g_i(B_j)}|\nabla_{H} u|^s\,dx
\biggl)^{1/s}
\\
\leq \Lambda K^{\nu-1}_{p,p}(f;D)N^{(s-1)/s}(f,D)
\biggr(\int\limits_{f^{-1}(B_j)}|\nabla_{H} u|^s\,dx
\biggl)^{1/s}.
\nonumber
\end{multline}
If $q<p$ then applying H\"older inequality twice we have
\begin{multline}
\biggr( \int\limits_{B_j}|\nabla_{H} v|^s\,dz
\biggl)^{1/s}
\leq \Lambda l^{(s-1)/s} \biggl[\sum\limits_{i=1}^{l}
\biggl(\int\limits_{g_i(B_j)}K_p(x;f)^{\frac{pq}{p-q}}~dx
\biggr)^{\frac{r-s}{r}} \cdot
\biggr(\int\limits_{g_i(B_j)}|\nabla_{H} u|^r\,dx\biggl)^{s/r}
\biggr]^{1/s}
\\
\leq \Lambda l^{(s-1)/s} \biggl[\biggl(\sum\limits_{i=1}^{l}
\int\limits_{g_i(B_j)}K_p(x;f)^{\frac{pq}{p-q}}~dx
\biggr)^{\frac{r-s}{r}} \cdot
\biggl(\sum\limits_{i=1}^{l}\int\limits_{g_i(B_j)}|\nabla_{H}
u|^r\,dx \biggr)^{s/r}
\biggr]^{1/s}\\
\leq \Lambda
l^{(s-1)/s}\|K_p(x;f)|L_{\frac{pq}{p-q}}(f^{-1}(B_j))\|^{\nu-1}
\biggr(\int\limits_{f^{-1}(B_j)}|\nabla_{H} u|^r\,dx
\biggl)^{1/r},
\nonumber
\end{multline}
where it was taken into account the equality
$\frac{rs}{r-s}=\frac{pq}{(p-q)(\nu-1)}$.

Since $l\leq N(f,D)$, we obtain 
\begin{multline}
\int\limits_{f(D)}|\nabla_{H}v|^s\,dz \leq
\sum\limits_{j=1}^{\infty}\int\limits_{B_j}|\nabla_{H} v|^s\,dx\\
\leq\sum\limits_{j=1}^{\infty}\Lambda^s
(N(f,D))^{(s-1)}\|K_p(x;f)|L_{\frac{pq}{p-q}}(f^{-1}(B_j))\|^{s(\nu-1)}
\biggr(\int\limits_{f^{-1}(B_j)}|\nabla_{H} u|^s\,dx
\biggl)^{s/r}.
\nonumber
\end{multline}

Applying the H\"older inequality we have
\begin{multline}
\int\limits_{f(D)}|\nabla_{H} v|^s\,dz \\
\leq \Lambda^s\bigl(N(f,D)\bigr)^{s-1}
\biggl(\sum\limits_{j=1}^{\infty}\|K_p(x;f)|L_{\frac{pq}{p-q}}(f^{-1}(B_j))\|^{\frac{rs}{r-s}(\nu-1)}\biggl)^{\frac{r-s}{r}}\biggl(\sum\limits_{j=1}^{\infty}\int\limits_{f^{-1}(B_j)}|\nabla_{H}
u|^r\,dx\biggr)^{\frac{s}{r}}\\
=\Lambda^s\bigl(N(f,D)\bigr)^{s-1}
\biggl(\sum\limits_{j=1}^{\infty}\|K_p(x;f)|L_{\frac{pq}{p-q}}(f^{-1}(B_j))\|^{\frac{pq}{p-q}}\biggl)^{\frac{r-s}{r}}
\biggl(\sum\limits_{j=1}^{\infty}\int\limits_{f^{-1}(B_j)}|\nabla_{H}
u|^r\,dx\biggr)^{\frac{s}{r}}
\\
\leq\Lambda^s\bigl(N(f,D)\bigr)^{s-1}
\|K_p(x;f)|L_{\frac{pq}{p-q}}\bigr(\bigcup\limits_{j=1}^{\infty}f^{-1}(B_j)\bigl)\|^{\frac{pq}{p-q}\frac{r-s}{r}}
\biggl(\int\limits_{\bigcup\limits_{j=1}^{\infty}f^{-1}(B_j)}|\nabla_{H}
u|^r\,dx\biggr)^{\frac{s}{r}}\\
\leq\Lambda^s \bigl(N(f,D)\bigr)^{s-1}\bigl(K_{p,q}(f;D)\bigr)^{s(\nu-1)}
\biggl(\int\limits_{D}|\nabla_{H} u|^r\,dx\biggr)^{\frac{s}{r}}.
\nonumber
\end{multline}
Theorem~3 is proved.

\vskip 0.3 cm

{\bf Corollary~1.} Let $f:\Omega\to\mathbb G$ be a
mapping with bounded $(p,q)$-distortion, $\nu-1<q\leq p<\infty$.
If $E=(A,C)$ is a condenser in the domain $\Omega$ such that
$\overline A\subset\Omega$, $C$  is a compact in $A$ and
$N(f,A)<\infty$, then
$$
\bigl(\cp_sf(E)\bigr)^{1/s}\leq\frac{(K_{p,q}(f;\Omega))^{\nu-1}(N(f,A))^{(s-1)/s}}{M(f,C)}\bigl(\cp_r
E)\bigr)^{1/r},
$$
where $r=\frac{q}{q-(\nu-1)}$ and $s=\frac{p}{p-(\nu-1)}$.

\vskip 0.3 cm

{\sc Proof.} Put in definition (2) of the function $v$ the
value $\Lambda=M^{-1}(f,C)$, where
$$
M(f,C)=\inf\limits_{x\in f(C)}\sum\limits_{z\in f^{-1}(x)\cap
C}i(z,f).
$$
Since $u$ is an admissible function for condenser $E(A,C)$, then
$v(x)\geq 1$ in the points $x\in f(C)$. Indeed, let $x\in f(C)$
and $f^{-1}(x)\cap C=\{z_1,z_2,\ldots,z_k\}$. Then
\begin{multline}
v(x)=\frac{1}{M(f,C)}\sum\limits_{z\in f^{-1}(x)}i(z,f)u(z)\\
\geq
\frac{1}{M(f,C)}\sum\limits_{l=1}^ki(z_l,f)u(z_l)
\geq
\frac{1}{M(f,C)}\sum\limits_{l=1}^ki(z_l,f) \geq 1
\nonumber
\end{multline}
since $u(z_l)\geq 1$.  Properties of the
function $v(x)$, proved in Theorem~3, imply that
$v$ is an admissible function for the condenser
$f(E)=(f(A),f(C))$. From here we have
\begin{multline}
\bigl(\cp_sf(E)\bigr)^{1/s} \leq
\biggl(\int\limits_{f(A)}|\nabla_{H}v|^s\,dx\biggr)^{1/s}\\
\leq
\frac{(K_{p,q}(f;\Omega))^{\nu-1}(N(f,A))^{(s-1)/s}}{M(f,C)}
\biggl(\int\limits_{A}|\nabla_{H} u|^r\,dx\biggr)^{1/r}.
\nonumber
\end{multline}
Since $u$ is an arbitrary function Corollary~1 is proved.

\vskip 0.3 cm

{\bf Corollary~2.} Let $f:\Omega\to\mathbb G$ be a
mapping with bounded $(p,q)$-distortion, $\nu-1<q\leq p<\infty$,
and $A=U(z,f,r_0)$ be a normal neighborhood of the point $z\in
\Omega$. Then for the condenser $E=(U(z,f,r_0),\overline U(z,f,r))$,
$0<r<r_0$, the estimate
$$
\bigl(\cp_s
f(E)\bigr)^{1/s}\leq\frac{K_{p,q}(f;\Omega)^{\nu-1}}{i(z,f)^{1/s}}
\bigl(\cp_r E\bigr)^{1/r}
$$
holds where $r=\frac{q}{q-(\nu-1)}$ and $s=\frac{p}{p-(\nu-1)}$.

\vskip 0.3 cm

{\sc Proof.} We note that $E=(U(z,f,r_0),\overline U(z,f,r))$ is
a normal condenser satisfying the conditions of Lemma~2. Hence
$N(f,A)=i(z,f)=\mu(f,A)=M(f,C)$. The desired inequality follows
from Corollary~1.

\vskip 0.3 cm

{\bf Remark~2.} In the  Euclidean space $\mathbb R^n$ and
$p=q=n$ Corollaries $1$ and $2$ were established in [27].

\vskip 0.3 cm

{\bf Proposition~2.} Let $f:\Omega\to\mathbb G$ be a
mapping with bounded $(p,q)$-distortion, $\nu-1<q\leq p<\infty$.
If $E=(A,C)$ is a condenser in the domain $\Omega$ such that
$\overline A\subset\Omega$,  $A$ is bounded, and $C$ is a
compact set in $A$, then
\begin{equation}
\bigl(\cp_sf(E)\bigr)^{1/s}\leq
K_{p,q}(f;\Omega)^{\nu-1}\bigl(\cp_r E\bigr)^{1/r}
\end{equation}
where $r=\frac{q}{q-(\nu-1)}$ and $s=\frac{p}{p-(\nu-1)}$.

\vskip 0.3 cm

{\sc Proof.} For $(p,q)$-quasiconformal mappings the inequality
(6) was proved in [55]. In the case $\mathbb G=\mathbb R^n$ the
inequality was proved in [48]. We describe the basic steps
of the proof in our situation. Since the closure $\overline A$ is
a compact set then $N(f,A)<\infty$. By a non-negative function
$u\in W^1_{\infty}(A)\cap C_0(A)$, we define the push-forward
function
$$
f_{\ast}u(z)=
\begin{cases} \sup\limits_{y\in f^{-1}(z)}u(y),\quad& z\in
f(\supp u),
\\
0, &z\notin f(\supp u).
\end{cases}
$$
 This operator
has the following properties, proofs of which are based on
Theorem~3:

1) the function $f_{\ast}u$ is continuous and $\supp
f_{\ast}u\subset f(\supp u)$;

2) $f_{\ast}:W^1_{\infty}(A)\cap C_0(A)\to W^1_{s}(\mathbb G)\cap
C_0(\mathbb G)$;

3) $\Bigl(\int\limits_{\mathbb
G}|\nabla_{H}f_{\ast}u|^s\,dx\Bigr)^{1/s}\leq
K_{p,q}(f;\Omega)^{\nu-1} \Bigl(\int\limits_{\mathbb
G}|\nabla_{H}u|^r\,dz\Bigr)^{1/r}$;

4) if a function $u$ is an admissible for the condenser $E=(A,C)$,
then  $f_{\ast}u$ is an admissible function for the condenser
$f(E)=(f(A),f(C))$.

From the last two properties we can to obtain the inequality (6).

\vskip 0.3 cm

The next theorem, giving an estimate for the local distortion of
the metrics under the mappings with bounded $(p,q)$-distortion
was proved in [42] for the Euclidean space $\mathbb R^n$. For
mappings of Carnot groups, the proof remain the same, taking into
account the lower estimate for the $q$-capacity, $\nu-1<q<\nu$
[55].

\vskip 0.3 cm

{\bf Theorem~4.} Let $f:\Omega\to\mathbb G$ be a  mapping with
bounded $(p,q)$-distortion where $\nu-1<q<\nu$. Consider an
arbitrary point $x\in\Omega$ and a neighborhood
$U(x,f,t_0),\,t_0=\min\{t_x,1\}$. Then, for every point $y\in
U(x,f,t)$, $t\leq t_0^4$ at $p=\nu$, and $t\leq t_0/4$ at
$p<\nu$, the following inequalities
$$
d(x,y)^{\frac{\nu-q}{q}}\leq 
\begin{cases}
C^{\frac{1}{\nu}}(\ln\frac{1}{d(f(x),f(y))})^{\frac{1-\nu}{\nu}}
\quad&\text{at $p=\nu$},\\
C^{\frac{1}{p}}(d(f(x),f(y))^{\frac{\nu-p}{p}} \quad&\text{at
$p<\nu$},
\end{cases}
$$
hold  where $C=K_{p,q}(f;U(x,f,t_0))\sup\limits_{y\in \mathbb
G}N(y,f,U(x,f,t_0))$.

\vskip 0.3 cm

\bigskip

\centerline{\bf 3.~Liouville type theorem and removable sets}

\bigskip

Let $K\subset \Omega$ be a compact set. We define the capacity of the
compact $\cp\bigl(K;W^1_p(\Omega)\bigr)$ in the space
$W^1_p(\Omega)$ as
$$\inf\{\|f\mid L_p(\Omega)\|^p+
\|\nabla_{H}f\mid L_p(\Omega)\|^p: \ f \in C_{0}^{\infty}(\Omega)\
\text{and}\ f\geq 1 \ \text{on}\ K\}.
$$
The capacity, defined initially on compact sets, extends by a
standard way on arbitrary sets (see, for example, [49,\,50],
where properties of the capacity are established).

\vskip 0.3 cm

{\bf Theorem~5.} Let $f:\mathbb G\to\mathbb G$ be a mapping with
bounded $(p,q)$-distortion, $\nu-1<q\leq p\leq\nu$. Then
$\cp\bigl(\mathbb G\setminus f(\mathbb G);W^1_s(\mathbb G)\bigr)=0$ where
$s={\frac{p}{p-(\nu-1)}}$.

\vskip 0.3 cm

{\sc Proof.} Indeed, fix a compact set  $C$ in $\mathbb G$ with a
nonempty interior, and  a sequence of open bounded sets $A_k\supset
C$  exhausting $\mathbb G$. For the condenser $E_k=(A_k,C)$, the
estimate (6)
$$
\bigl(\cp_sf(E_k)\bigr)^{1/s}\leq K_{p,q}^{\nu-1}(f;\mathbb
G)\bigl(\cp_r E_k\bigr)^{1/r}
$$
holds. Since  the right-hand side goes to zero as $k\to\infty$,
then the left-hand side does the same. Note that
$$
\lim\limits_{k\to\infty}\cp_sf(E_k)=\cp_s(f(\mathbb G),f(C);\mathbb
G)=0.
$$
Hence
$$
\cp\bigl(\mathbb G\setminus f(\mathbb G);W^1_s(\mathbb G)\bigr)\leq
\cp_s(\mathbb G\setminus f(\mathbb G),f(C);\mathbb G)=\cp_s(f(\mathbb
G),f(C);\mathbb G)=0.
$$

From here it is obviously follows, that $\cp\bigl(\mathbb G\setminus
f(\mathbb G);W^1_s(\mathbb G)\bigr)=0$.

\vskip 0.3 cm

{\bf Corollary~3.} Let $f:\mathbb G\to\mathbb G$ be a mapping with
bounded distortion, and $\cp\bigl(\mathbb G\setminus f(\mathbb
G);W^1_{\nu}(\mathbb G)\bigr)>0$. Then $f$ is a constant mapping.

\vskip 0.3 cm

{\bf Corollary~4.} Let $f:\mathbb G\to\mathbb G$ be a mapping with
bounded $(p,q)$-distortion, $p,q\in(\nu-1,\nu)$. Then $f(\mathbb
G)=\mathbb G$.

\vskip 0.3 cm

In order to formulate the next assertion, we recall that the
space $w^1_p(\mathbb G)$, $1<p<\nu$, is defined as the completion of
the space $C_{0}^{\infty}(\mathbb G)$ with respect to the  norm
$L^1_p(\mathbb G)$. The capacity of the compact $K$ in the space
$w^1_p(\mathbb G)$ is defined by the following way:
$$
\cp\bigl(K;w^1_p(\mathbb G)\bigr)= \inf\{\|\nabla_{H}f\mid L_p(\mathbb
G)\|^p: \ f \in C_{0}^{\infty}(\mathbb G)\ \text{and}\ f\geq 1 \
\text{on}\ K\}.
$$
It is proved in [49,\,50], that the families of the sets of
capacity zero in the spaces $w^1_p(\mathbb G)$ and $W^1_p(\mathbb G)$
coincide.

Let $E\subset\mathbb G$ be a closed set, $\nu$-capacity of which is
positive: $\cp\bigl(E;$ $W_{\nu}^1 (\mathbb G)\bigr)>0$. We say, that
set $E$ has an essentially positive capacity at the point $x\in
E$ if $\cp\bigl(E\cap B(x,r);W_{\nu}^1(B(x,2r))\bigr)>0$ for
every $r\in(0,1)$. It is easy to show that if
$\cp\bigl(E;W_{\nu}^1(\mathbb G)\bigr)>0$ then collection $\widetilde
E$ of the points of the essentially positive capacity for the set
$E$ is a nonempty set. We note that the set $\widetilde E$ is
closed. Hence, there exists point $x_0\in \widetilde E$ the closest to zero:
$d(x_0)=\inf\{d(x):x\in\widetilde E\}$.
Consider the intersection $E_0=\widetilde E\cap
B(x_0,1)$. The capacity $\cp\bigl(E_0;W_{\nu}^1(\mathbb G)\bigr)$ is
positive.

The following two lemmas were established in [53].

\vskip 0.3 cm

{\bf Lemma~11.} Let $C$ be a continuum in $\mathbb G$ and $\diam
C\geq\alpha>0$.

$1)$ The case $\nu-1<p<\nu$. For every $\alpha>0$, there exists $\delta>0$
such that if $\diam C\geq\alpha>0$ then $\cp\bigl(C;w_p^1(\mathbb
G)\bigr)>\delta$.

$2)$ The case $p=\nu$. If $E\subset\mathbb G$ is a compact set and its capacity
$\cp\bigl(E;W_{\nu}^1(\mathbb G)\bigr)$ is positive, then for every
$\alpha>0$ and $d>0$, there exists $\delta>0$ such that
$\cp_{\nu}(\mathbb C E, C)>\delta$ under conditions $\diam
C\geq\alpha>0$ and $\dist(C,E_0)\leq d$.

\vskip 0.3 cm

As a corollary, we obtain the next assertion.

\vskip 0.3 cm

{\bf Lemma~12.} Let $\nu-1<p<\nu$, $A$ is a bounded open set
in $\mathbb G$, and $C$ is a continuum in $A$. Then for every
$\alpha>0$ there exists a number $\varepsilon>0$ such that if
$\diam C\geq\alpha$ then $\cp_p\bigl(A,C\bigr)\geq\varepsilon$.

\vskip 0.3 cm

{\sc Proof.} It is enough to apply the inequality
$\cp\bigl(C;w^1_p(\mathbb G)\bigr)\leq\cp_p(A,C)$.

\vskip 0.3 cm

{\sc Corollary~5.} Consider the family $f:\Omega\to\mathbb G$,
$\Omega\subset\mathbb G$, of mappings with bounded
$(p,q)$-distortion and with a coefficient of distortion not
greater than a fix number $K$.

$1)$ If $p\geq q>\nu$, then the given family of the mappings is
locally uniformly continuous.

$2)$ If $p=q=\nu$ and the family of the mappings $f:\Omega\to\mathbb
G\setminus E$, where $E$ is a closed set of positive capacity
$\cp\bigl(E;W_{\nu}^1(\mathbb G)\bigr)$, is locally uniformly
bounded, then it is locally equicontinuous.

\vskip 0.3 cm

{\sc Proof.} Fix a compact set $A\subset\Omega$ and an open
bounded set $U\supset A$ with $\overline U\subset\Omega$. There
exists $R_0$ such that capacity $\cp_r(U,\overline{B(x,R)})$ is
uniformly small with respect to all $R\in(0,R_0)$ and $x\in A$,
$r=\frac{q}{q-(\nu-1)}$. By inequality (6), we have the estimate
$$
\bigl(\cp_{s}(f(U),f(\overline {{B(x,r)}}))\bigr)^{1/s} \leq
K_{p,q}(f;\Omega)^{\nu-1} \bigl(\cp_{r}(U,\overline
{B(x,r)})\bigr)^{1/r}, \  s={p/(p-\nu+1)}.
$$
Now the required assertions follow from Lemma~11 and Lemma~12.

\vskip 0.3 cm

{\bf Theorem~6.} Let $f:\Omega\setminus F\to\mathbb G$ be a
 mapping with bounded $(p,q)$-distortion, $p\geq
q\geq\nu$, and $F$ be a closed set in the domain $\Omega$,
$\cp\bigl(F;W^1_s(\mathbb G)\bigr)=0$, $s=p/(p-(\nu-1))$. Then

$1)$ in the case $p\geq q>\nu$: the mapping $f$ extends to a continuous
mapping $\tilde f:\Omega\to\mathbb G$;

$2)$ in the case  $p=q=\nu$: if $\cp\bigl(\mathbb Cf(\Omega\setminus
F);W^1_{\nu}(\mathbb G)\bigr)>0$ then the mapping $f$ extends to a
continuous mapping $\tilde f:\Omega\to\mathbb G\cup{\infty}$.
$($Hereafter  $\mathbb G\cup\{\infty\}$ is the one-point
compactification $\mathbb G$ with the standard topology.$)$

\vskip 0.3 cm

{\sc Proof.} {\sl The case} $p\geq q>\nu$. Let $x_0\in F$ and the
mapping $f$ has no limit at the point $x_0$. Then in a some ball
$B(x_0,R)\subset \Omega$ exists two sequences of points $\{x_i\}$
and $\{x^{\prime}_{i}\}$, tends to the point $x_0$, but for every
index $i$ the inequality $d(f(x_i),f(x^{\prime}_i))\geq\alpha>0$
holds. Since $\cp_{p}(F)=0$, then the sets $B(x_0,r_i)\setminus F$
are connected, here $r_i=\max\{d(x_0,x_i),d(x_0,x^{\prime}_i)\}$.
The points $x_i$ and $x^{\prime}_i$ can be connected by a
horizontal curve $C_i$, lying in the set $B(x_0,r_i)\setminus F$.
Since $\diam(f(C_i))\geq\alpha$, then by Lemma~11,
$\cp_{s}(f(B(x_0,R)\setminus F),f(C_i))\geq\varepsilon^s>0$ for
every index $i$. Let open sets $A_j$ such, that them compact
embedded in $\Omega$, $A_j\subset B(x_0,R)\setminus F$ and as
$j\to\infty$ exhausts $B(x_0,R)\setminus F$. Then by inequality
(6) we have
$$
\bigl(\cp_{s}(f(A_j),f(C_i))\bigr)^{1/s} \leq
K_{p,q}(f;\Omega)^{\nu-1}\bigl(\cp_{r}(A_j,C_i)\bigr)^{1/r}.
$$
Passing to the limit as $j\to\infty$, we arrive to the relations
\begin{multline}
\varepsilon \leq\bigl(\cp_{s}(f(B(x_0,R)\setminus
F),f(C_i))\bigr)^{1/s} \leq K_{p,q}(f;\Omega)^{\nu-1}
\bigl(\cp_{r}(B(x_0,R)\setminus F,C_i)\bigr)^{1/r}
\\
\leq K_{p,q}(f;\Omega)^{\nu-1}
\bigl(\cp_{r}(B(x_0,R),C_i)\bigr)^{1/r}.
\end{multline}
Since $ \cp_{r}(B(x_0,R),C_i) \leq \theta\bigl({\frac{r-1}{r-\nu}}\bigl(
R^{\frac{r-\nu}{1-r}}-r_i^{\frac{r-\nu}{1-r}}\bigr)\bigr)^{1-r}, $
then the capacity in the right side of the inequalities (7) tends
to zero, as $i\to \infty$.

We obtained the contradiction, since the left side of the
inequality (7) is delimited from zero, but the right side tends
to zero, while $i\to \infty$.

{\sl The case} $p=q=\nu$, is divided in two subcases 
depending on whether 
$$
\min(\dist(f(x_i), E_0), \dist(f(x^{\prime}_i), E_0)\to \infty,
$$
or not (here $E=\mathbb Cf(\Omega\setminus F)$, and definition
of the set $E_0$ before Lemma~11). If it goes to $\infty$, then we can
assume that $\lim\limits_{x\to x_0}f(x)=\infty$. If it does not
go then the proof does not differ from the previous case.
Theorem~6 is proved.

\vskip 0.3 cm

{\bf Corollary~6.} Let $f:\Omega\to\mathbb G$ be a
mapping with bounded $(p,q)$-distortion, $p\geq q\geq\nu$, and a
point $b\in\partial \Omega$ is isolated.
Then

$1)$ in the case $p\geq q>\nu$: the mapping $f$ has a continuous
extension $\tilde f:\Omega\cup\{b\}\to\mathbb G$;

$2)$ in the case $p=q=\nu$: if $\cp\bigl(\mathbb Cf(\Omega);W^1_{\nu}(\mathbb
G)\bigr)>0$, then the mapping $f$ extends to a continuous mapping
$\tilde f:\Omega\cup\{b\}\to\mathbb G\cap\{\infty\}$.

\vskip 0.3 cm

Proof of  Corollary~6 is obvious.

{\bf Corollary~7.} Let $f:\Omega\to\mathbb G$ be a nonconstant
mapping with bounded distortion, a point $b\in\partial \Omega$ be
isolated, and the mapping $f$ have no limit at the point $b$.
Then $\cp\bigl(\mathbb Cf(U\setminus\{b\});W^1_{\nu}(\mathbb
G)\bigr)=0$ for every neighborhood $U\subset\Omega\cup\{b\}$ of
the point $b$. Furthermore, there exists ${\sigma}$-set
$E\subset\mathbb G\cup\{\infty\}$, having the capacity zero in the
space $W^1_{\nu}(\mathbb G)$, such that
$N(z,f,U\setminus\{b\})=\infty$ for every point $z\in (\mathbb
G\cup\{\infty\})\setminus E$ and for every above-mentioned
neighborhood of the point $b$.

\vskip 0.3 cm

{\sc Proof.} For coming to a contradiction, we assume
$$
\cp\bigl(\mathbb Cf(U\setminus\{b\});W^1_{\nu}(\mathbb G)\bigr)>0.
$$
By Corollary~6, the mapping $f$ has the limit at the point $b$. It
contradicts to a condition of the corollary. Indeed, let $k_0$ be
a number such that $1/k_0<\dist(b,\partial U)$. The set $\mathbb
Cf(B(b,1/k)\setminus\{b\})$ has the capacity zero in the space
$W^1_{\nu}(\mathbb G)$ for every $k\geq k_0$. Then the set
$$
E=\bigcup\limits_{k=k_0}^{\infty}\mathbb Cf(B(b,1/k)
$$
in view of subadditivity of the capacity, has also capacity zero
in the space $W^1_{\nu}(\mathbb G)$. Moreover, the preimage of every
point $z\in\mathbb CE$ belongs to the balls $B(b,1/k)$ for $k\geq
k_0$. Hence, there exists a sequence pairwise different points $x_i\in
B(b,1/k_i)$ such that $f(x_i)=z$ for every $i$. Corollary is
proved.

\vskip 0.3 cm

{\bf Theorem~7.} Let $f:\Omega\to\mathbb G$ be a
mapping with bounded $(p,q)$-distortion, $\nu-1<q\leq p<\infty$,
and a point $b\in\partial \Omega$ be isolated. Then, if the
mapping $f$ permits a continuous extension $\tilde
f:\Omega\cup\{b\}\to\mathbb G$, then $\tilde f$ is  a mapping with
bounded $(p,q)$-distortion.

\vskip 0.3 cm

{\sc Proof.}  It is obvious that the extension $\tilde f$ belongs
to the class $\ACL$,  satisfies the relation
$$
|D_H\tilde f(x)| \leq K_p(x;f)J(x,\tilde f)^{\frac{1}{p}}
$$
almost everywhere, and has Luzin condition $\mathcal N$. We have to show,
that the function $|D_H\tilde f|^{q}$ is integrable in the
neighborhood of the point $b$. Since the mapping $f$ is
discrete, then there exists a ball $B(b,r)\Subset\Omega$ such that its
image is bounded and $\partial B(b,r)\cap \tilde f^{-1}(\tilde
f(b))=\emptyset$. Let $U$ be $f(b)$-component of the set $\mathbb
C\tilde f(\partial B(b,r))$, and $V$ be $b$-component of the set
$\tilde f^{-1}(U)$. By change of variable formula (see Proposition 1), we have
\begin{multline} \int\limits_{V}|D_{H}\tilde f(x)|^{q}\,dx \leq
\int\limits_{V}(K_p^q(x;f)J(x,\tilde f))^{q/p}\,dx\\
\leq
\biggl(\int\limits_{V}(K_p(x;f)^{\frac{pq}{p-q}}\,dx\biggr)^{\frac{p-q}{p}}
\biggl(\int\limits_{V}J(x,\tilde f))\,dx\biggr)^{q/p}\\
\leq K_{p,q}(f;V)^{q}\biggl(\int\limits_{U}N(z,\tilde
f,V)\,dz\biggr)^{q/p} < \infty.
\nonumber
\end{multline}
Theorem~7 is proved.

\vskip 0.3 cm

\bigskip

\centerline{\bf 4. Geometrical definition of $(p,q)$-quasiregular mappings}

\bigskip

Boundedness of a linear distortion is an equivalent geometric
characteristic for mappings with bounded distortion
($\nu$-distortion).

The following description of mappings in the case $\mathbb G=\mathbb
R^n$ was formulated in [27]. A nonconstant mapping
$f:\Omega\to\mathbb G$, $\Omega\subset\mathbb G$, is called
quasiregular, if

{(1)} $f$ is continuous, open, discrete, and sense-preserving
at the points of a local homeomorphism;

{(2)} The value $H(x,f)=\overline{\lim\limits_{r\to
0}}\frac {\max\limits_{d(x,y)=r}d(f(x),f(y))}
{\min\limits_{d(x,y)=r}d(f(x),f(y))}$ is local bounded in
$\Omega$;

{(3)}  For all points $x\in\Omega\setminus B_f$, where $B_f$ is a
branch set of $f$, the relation $H(x,f)\leq K$ holds where $K$ is independent
of a point $x$.

Analytic properties of quasiregular mappings are formulated in
the next assertion.

\vskip 0.3 cm

{\bf Theorem~8 \rm[53]}. Every nonconstant quasiregular
mappings $f:\Omega\to\mathbb G$, defined on domain
$\Omega\subset\mathbb G$, is a mapping with bounded distortion.

\vskip 0.3 cm

{\sc Proof.} $\ACL$-property of the mapping $f$ was proved in
[52]. The integrability of the horizontal differential
follows from the change of variable formula [45].

\vskip 0.3 cm

In the similar way we can to describe geometrically the
mappings with bounded $(p,q)$-distortion.

A nonconstant mapping  $f\:\Omega\to\mathbb G$, $\Omega\subset\mathbb
G$, is called $(p,q)$-quasiregular if

{(1)} $f$ is continuous open discrete and sense-preserving
at the points of local homeomorphism;

{(2)} in the domain $\Omega$, there exists a bounded quasiadditive set function $\Phi$ defined on open bounded subsets of $\Omega$ such that the value
$$ H_{p,q}(x,f)={\varlimsup\limits_{r\to 0}}
\frac{L_{f}^{p}(x,r)r^{{\nu}-p}} {|f(B(x,\lambda r))|}
\bigg/\biggl(\frac{\Phi(B(x,\lambda r))}{|B(x,r)|}\biggr)^{\frac{p-q}{q}}
$$
is locally bounded in $\Omega$. Here $\lambda>1$ is a fixed
number, and $L_{\varphi}(x,r)=\max\{d(f(y),f(x)):
d(y,x)=r\}$ by the condition that $r<\dist(x,\partial D)$;

{(3)}  For all points $x\in\Omega\setminus B_f$, where $B_f$ is the
branch set of $f$, the relation $H_{p,q}(x,f)\leq K$ holds,
where $K$ is independent of  $x$;

{(4)} under $p,q\in(\nu-1,\nu)$  the mapping $f$ has Luzin
condition  $\mathcal N$.

\vskip 0.3 cm

{\bf Theorem~9.} Every  $(p,q)$-quasiregular
mapping $f :\Omega\to\mathbb G$, defined on the domain
$\Omega\subset\mathbb G$, is a mapping with bounded
$(p,q)$-distortion.

\vskip 0.3 cm

{\sc Proof.} 
Let a point $x\in D\setminus B_f$, then there exists a neighborhood $W$ of the point $x$ such that the restriction
$f : W\to f(W)$ be a homeomorphism. Then $f$ is $\ACL$-mapping on $W$ and is differentiable almost everywhere in $W$ (see, for details, [55]). Therefore $f$ belongs to the class $\ACL(\Omega)$ and is differentiable almost everywhere in $\Omega$ . The integrability of the local $p$-distortion follows now from the boundedness of the value $H_{p,q}(x,f)$.
Indeed, since $f$ is an $\ACL$ mapping differentiable a. e. in $\Omega$, we have
$$
H_{p,q}(x,f)=\frac{|Df|^p}{\lambda^\nu J(x,f)}\bigg/\biggl(\lambda^\nu\Phi'(x)\biggr)^{\frac{p-q}{q}}\leq K \quad\text{a.~e.~ in}\quad\Omega.
$$
Hence
$$
\frac{|Df|^p}{\lambda^\nu J(x,f)}\leq K \biggl(\lambda^\nu\Phi'(x)\biggr)^{\frac{p-q}{q}}\quad\text{a.~e.~ in}\quad\Omega.
$$
Since the volume derivative $\Phi'(x)$ integrable in $\Omega$ [57], we see that the value
$$
I_{p,q}(x,f)=\biggl(\frac{|Df|^p}{J(x,f)}\biggr)^{\frac{q}{p-q}}
$$
is integrable in $\Omega$ and $f$ is the mapping with bounded $(p,q)$ distortion.
\vskip 0.5cm

\centerline{REFERENCES}

\begin{enumerate}

\item
Astala K., Iwaniec T., Koskela P., Martin G. {\it Mappings of BMO-bounded distortion}//
Math. Ann. -- 2000. -- V.~317. -- P.~703--726.

\item
Ball J.~M. {\it Convexity conditions and existence theorems in nonlinear elasticity}//
Arch. Rat. Mech. Anal. -- 1976. -- V.~63. -- P.~337--403.

\item
Ball J.~M. {\it Global invertability of Sobolev functions and the interpretation of matter}//
Proc. Roy. Soc. Edinburgh -- 1981. -- V.~88A. -- P.~315--328.

\item
Dairbekov N.~S. {\it Mappings with bounded distortion on Heisenberg groups}//
Siberian Math. J. -- 2000. -- V.~41. -- N.~3. -- P.~465--486.

\item
Dairbekov N.~S. {\it Mapping with bounded distortion of two-step Carnot groups}//
Proceedings on Analysis and Geometry --   Novosibirsk: Sobolev Institute Press -- 2000. -- P.~122--155.

\item
Federer H. {\it Geometric measure theory} -- Berlin: Springer-Verlag. -- 1969.

\item
Folland G.~B and Stein E.~M. {\it Hardy spaces on homogeneous group} -- Princeton: Princeton Univ. Press. -- 1982.

\item
Gol'dshtein V., Gurov L. {\it Applications of change of variables operators for exact embedding theorems}//
Integral Equations Operator Theory -- 1994. -- V.~19. -- N.~1. -- P.~1-24.

\item
Gol'dshtein V., Ramm A.~G. {\it Compactness of the embedding operators for rough domains}//
Math. Inequalities and Applications -- 2001. -- V.~4. -- N.~1. -- P.~127-141.

\item
Gol'dshtein V., Ukhlov A. {\it Weighted Sobolev spaces and embedding theorems}//
Trans. Amer. Math. Soc. (to appear)

\item
Gromov M. {\it Carnot--Caratheodory spaces seen from within}//
Sub-Reimannian geometry -- 1996. -- Basel: Birkh\"auser. -- P.~79--323.

\item
Heinonen~J., Holopainen~I. {\it Quasiregular maps on Carnot groups}//
Geom. Anal. -- 1997. -- V.~7. -- N.~1. -- P.~109--148.

\item
Heinonen J., Kilpel\"ainen T., Martio O. {\it Nonlinear potential theory of degenerate elliptic equations} --
Oxford--New-York--Tokio: Clarendon Press. 1993.

\item
Heinonen J., Koskela P. {\it Sobolev mappings with integrable dilatation}//
Arch. Rat. Mech. Anal. -- 1993. -- V.~125. -- P.~81--97.

\item
Hencl S., Koskela P. {\it Mappings of finite distortion: openness and discteteness for quasilight mappings}//
Ann. Inst. H. Poincare Anal. Non Lineaire -- 2005. -- V.~22. -- P.~331--342.

\item
Iwanec T., Koskela P., Martin G. {\it Mappings of finite distortion: Monotonicity and continuity}//
Invent. Math. -- 2001. -- V.~114. -- N.~3. -- P.~507--531.

\item
Iwanec T., Koskela P., Martin G. {\it Mappings of BMO-distortion and Beltrami type operators}//
J. Anal. Math. -- 2002. -- V.~88. -- P.~337--381.

\item
Iwanec T., \v{S}ver\'{a}k V. {\it On mappings with integrable dilatation}//
Proc. Amer. Math. Soc. -- 1993. -- V.~118. -- P.~181--188.

\item
Kauhanen J., Koskela P., Mal\'{y} J. {\it Mappings of finite distortion: Discteteness and openness}//
Arch. Rat. Mech. Anal. -- 2001. -- V.~160. -- N.~2. -- P.~135--151.

\item
Kauhanen J., Koskela P., Mal\'{y} J. {\it Mappings of finite distortion: Condition N}//
Michigan Math. J. -- 2001. -- V.~49. -- N.~1. -- P.~169--181.

\item
Kauhanen J., Koskela P., Mal\'{y} J. {\it Mappings of finite distortion: Sharp Orlicz-conditions}
Rev. Mat. Iberoamericana -- 2003. -- V.~19. -- N.~3. -- P.~857--872.

\item
Koranyi A., Reimann H.~M. {\it Foundations for the theory of quasiconformal mappings on the Heisenberg group}//
Adv. Math. -- 1995. -- V.~111. -- N.~1. -- P.~1--87.

\item
Koskela P., Onninen J. {\it Mappings of finite distortion: capacity and modulus inequalities}//
J. Reine Angew. Math. -- 2006. -- V.599. -- P.~1--26.

\item
Manfredi J. {\it Weakly monotone functions}// J. Geom. Anal. -- 1994. -- V.~4. -- P.~393--402.

\item
Manfredi J., Villamor E. {\it Mappings with integrable dilatation in higher dimensions}//
Bull. Amer. Math. Soc. -- 1995. -- V.~32. -- N.~2. -- P.~235--240.

\item
Markina I.~G., Vodop'yanov S.~K. {\it On value distribution for quasimeromorphic mappings on $\mathbb H$-type Carnot groups}//
Bulletin des Sciences Math\'ematiques. -- 2006. -- V.~130. -- P.~467--523.

\item 
Martio O. {\it A capacity inequality for quasiregular mappings}//
Ann. Acad. Sci. Fenn. Ser. A I Math. -- 1971. -- V.~474. -- P.~1--18.

\item
Martio O., Rickman S.,V\"ais\"al\"a J. {\it Definitions for quasiregular mappings}//
Ann. Acad. Sci. Fenn. Ser. A I Math. -- 1971. -- V.~448. -- P.~1--40.

\item
Mostow G.~D. {\it Strong rigity of locally symmetric spaces} -- Princeton: Princeton Univ. Press. 1973.

\item
Ottazzi A. {\it A sufficient condition for nonrigidity of Carnot groups}//
Mathematische Zeitschrift (to appear)

\item
Pansu P. {\it M\'etriques de Carnot--Carath\'eodory et quasiisom\'etries des espaces sym\'et\-ri\-qu\-es de rang un}//
Ann. Math. -- 1989. -- V.~129. -- P.~1--60.

\item
Reshetnyak Yu.~G. {\it Space mappings with bounded distortion} -- Novosibirsk: Nauka. 1982.

\item
Reshetnyak Yu.~G. {\it Sobolev-type classes of functions with values in a metric space}//
Siberian Math. J. -- 1997. -- V.~38. -- P.~657--675.

\item
Rickman S. {\it Quasiregular mapping} -- Berlin--Heidelberg--New-York: Springer - Verlag. 1993.

\item
Ryazanov V., Srebro U., Yakubov E. {\it BMO-quasiconformal mappings}//
J. Anal. Math. -- 2001. -- V.~83. -- P.~1--20.

\item
Ryazanov V., Srebro U., Yakubov E. {\it The Beltrami equation and FMO functions}//
Complex analysis and dynamical systems II, Contemp. Math. (Amer. Math. Soc., Providence, RIJ. Anal. Math.)
-- 2005. P.~357--364.

\item
Troyanov M., Vodop'yanov S.~K. {\it Liouville type theorems for mappings with bounded (co)-distortion}//
Ann. Inst. Fourier (Grenoble) -- 2002. -- V.~52. -- V.~6. -- P.~1753--1784.

\item
Ukhlov A. {\it Differential and geometrical properties of Sobolev mappings.}//
Matem. Notes -- 2004. -- V.~75. -- P.~291-294.

\item
Ukhlov A., Vodop'yanov S.~K. {\it Mappings associated with weighted Sobolev spaces}//
Contemporary Mathematics -- 2008. -- V.~ -- P.~.

\item
Vodop'yanov S.~K. {\it Monotone functions and quasiconformal mappings on Carnot groups}//
Siberian Math. J. -- 1996. -- V.~37. -- N.~6. -- P.~1113--1136.

\item
Vodop'yanov S.~K. {\it Mappings with bounded distortion and with finite distortion on Carnot groups}//
Siberian Math. J. -- 1999. -- V.~40. --- N.~4. -- P.~764--804.

\item
Vodop'yanov S.~K. {\it Topological and geometric properties of mappings with an integrable Jacobian in Sobolev classes. I.}// Siberian Math. J. -- 2000. -- V.~41. -- N.~4. -- P.~19--39.

\item
Vodop'yanov S.~K. {\it $\cal P$-Differentiability on Carnot Groups in Different Topologies and Related Topics}//
Proceedings on Analysis and Geometry -- 2000. -- Novosibirsk: Sobolev Institute Press. -- P.~603--670.

\item
Vodop'yanov S.~K. {\it Closure of classes of mappings with bounded distortion on Carnot groups}//
Mat. tr. -- 2002. -- V.~5. -- N.~2. -- P.~92--137.

\item
Vodop'yanov S.~K. {\it Differentiability of maps of Carnot groups of Sobolev classes}//
Sbornik: Mathematics -- 2003. -- V.~194. -- N.~6. -- P.~857--877.

\item
Vodop'yanov S.~K. {\it Composition Operators on Sobolev Spaces}//
Contemporary Ma\-thematics -- 2005. -- V.~382. -- P.~327-342.

\item
Vodopyanov~S.~K. {\it Foundations of the Theory of Mappings with Bounded Distortion on Carnot Groups}//
Contemp. Math. -- 2007. --  V.~424. -- P.~303--344.

\item
Vodop'yanov S.~K. {\it Topological and geometrical properties of Sobolev classes with integrable Jacobian. II.}//
Sibirsk. Mat. Zh. (to appear)

\item
Vodop'yanov S.~K., Chernikov V.~M. {\it Sobolev spaces and hypoelliptic equations}//
Proceeding Institute of Mathematics -- Siberian Branch of RAS, Novosibirsk -- 1995. -- 3--64.

\item
Vodop'yanov S.~K., Chernikov V.~M. {\it Sobolev spaces and hypoelliptic equations. II}//
Siberian Advances in Mathematics -- 1996. -- N.~4. -- P.~64--96.

\item
Vodop'yanov S.~K., Gol'dshtein V.~M. {\it Quasiconformal mappings and spaces of functions with generalized first 
derivatives}// Siberian Math. J. -- 1976. -- V.~17. -- N.~3. -- P.~515--531.

\item
Vodop'yanov S.~K., Greshnov A.~V. {\it Analytical properties of quasiconformal mappings on Carnot groups}//
Siberian Math. J. -- 1995. -- V.~36. -- N.~4. -- P.~1317--1327.

\item
Vodop'yanov S.~K., Markina I.~G. {\it Local estimates for mappings with bounded $s$-distortion on Carnot groups}//
Proceeding of 12 Siberian School "Algebra, Geometry, Analysis and Mathematical Physics" -- Novosibirsk: Siberian Branch of RAS -- 1999. -- P.~27--43.

\item
Vodop'yanov S.~K., Ukhlov A.~D. {\it Approximately differentiable transformations and change of variables on nilpotent groups}// Siberian Math. J. -- 1996. -- V.~37. -- N.~1. -- P.~79--80.

\item
Vodop'yanov S.~K., Ukhlov A.~D. {\it Sobolev spaces and $(P,Q)$-quasiconformal mappings of Carnot groups}//
Siberian Math. J. -- 1998. -- V.~39. -- N.~4. -- P.~665-682.

\item
Vodop'yanov S.~K., Ukhlov A.~D. {\it Superposition operators in Sobolev spaces}//
Russian Mathematics (Iz. VUZ) -- 2002. -- V.~46. -- P.~11-33.

\item
Vodop'yanov S.~K., Ukhlov A. {\it Set functions and their applications in the theory of Lebesgue and Sobolev spaces. I}//
Siberian Advances in Mathematics. -- 2004. -- V.~14. -- P.~78-125.

\item
Vodop'yanov S.~K., Ukhlov A. {\it functions and their applications in the theory of Lebesgue and Sobolev spaces. II}//
Siberian Adv. Math. -- 2005. -- V.~15. -- N.~1. -- P.~91-125.

\item
Warhurst B. {\it Contact and quasiconformal mappings on real model filiform groups}//
Bull. Austral. Math. Soc. -- 2003. -- V.~68. -- P.~329--343.

\item
Warhurst B. {\it Jet spaces as nonrigid Carnot groups}//
J. Lie Theory  -- 2005. -- V.~15. -- P.~341--356.

\end{enumerate}

\end{document}